\documentclass{conm-p-l}

\usepackage{amssymb}
\usepackage[all]{xy}
\usepackage{times}
\usepackage{graphics}
\usepackage{array}
\usepackage{bbold}
\usepackage{bm}

\newtheorem{thm}{Theorem}
\newtheorem{prop}[thm]{Proposition}
\newtheorem{lem}[thm]{Lemma}

\theoremstyle{definition}

\newtheorem{example}[thm]{Example}
\newtheorem*{ack}{Acknowledgment}

\theoremstyle{remark}
\newtheorem*{notat}{Notations}
\newtheorem*{claim}{Claim}

\newtheorem*{remarks}{Remarks}

\DeclareMathOperator{\tr}{trace}
\DeclareMathOperator{\Tr}{Tr}
\DeclareMathOperator{\ch}{char}

\DeclareMathOperator{\diag}{diag}

\DeclareMathOperator{\GL}{GL}
\DeclareMathOperator{\PGL}{PGL}

\DeclareMathOperator{\M}{M}
\DeclareMathOperator{\N}{\mathcal{N}}

\DeclareMathOperator{\End}{End}

\DeclareMathOperator{\Aut}{Aut}
\DeclareMathOperator{\Gal}{Gal}
\DeclareMathOperator{\Id}{Id}

\newcommand{\Sy}{\mathcal{S}}

\DeclareMathOperator{\Or}{O}
\DeclareMathOperator{\Sp}{Sp}

\DeclareMathOperator{\Iso}{Iso}

\renewcommand{\k}{\mathbb{k}}

\newcommand{\into}{\subseteq}

\newcommand{\G}{\mathcal{G}}
\newcommand{\e}{\varepsilon}
\newcommand{\z}{\zeta}
\newcommand{\p}{\mathfrak{p}}
\newcommand{\bP}{\mathfrak{P}}

\renewcommand{\P}{\mathfrak{P}}

\newcommand{\cO}{\mathcal{O}}

\newcommand{\FF}{\mathbb{F}}
\newcommand{\FFq}{\mathbb{F}_{\!q}}
\newcommand{\CC}{\mathbb{C}}
\newcommand{\NN}{\mathbb{N}}
\newcommand{\ZZ}{\mathbb{Z}}
\newcommand{\QQ}{\mathbb{Q}}
\newcommand{\RR}{\mathbb{R}}
\bmdefine{\bmu}{\mu}

\newcommand{\wed}[1]{\textstyle\bigwedge^2{#1}}
\newcommand{\Sym}{\mathsf{S}}
\newcommand{\sym}{\Sym^2}

\begin{document}

\title[Finite matrix groups]{Orders of Finite Groups of Matrices}

\author{Robert M. Guralnick}
\address{Department of Mathematics, University of Southern California, Los Angeles, CA 90089-1113}
\email{guralnick@usc.edu}
\thanks{The first author was supported in part by NSF Grant DMS 0140578.}

\author{Martin Lorenz}
\address{Department of Mathematics, Temple University,
    Philadelphia, PA 19122-6094}
\email{lorenz@math.temple.edu}
\thanks{Research of the second author supported in part by a grant from the NSA}


\dedicatory{To Don Passman, on the occasion of his $65^\text{th}$
birthday}


\begin{abstract}
We present a new proof of a theorem of Schur's from 1905 determining
the least common multiple of the orders of all finite groups of
complex $n \times n$-matrices whose elements have traces in the
field $\QQ$ of rational numbers. The basic method of proof goes back
to Minkowski and proceeds by reduction to the case of finite fields.
For the most part, we work over an arbitrary number field rather
than $\QQ$. The first half of the article is expository and is
intended to be accessible to graduate students and advanced
undergraduates. It gives a self-contained treatment, following
Schur, over the field of rational numbers.
\end{abstract}

\maketitle


\section{Introduction}

\subsection{} 
How large can a finite group of complex $n \times n$-matrices be if
$n$ is fixed? Put differently: if $\G$ is a finite collection of
invertible $n \times n$-matrices over $\CC$ such that the product of
any two matrices in $\G$ again belongs to $\G$, is there a bound on
the possible cardinality $|\G|$, usually called the \emph{order} of
$\G$? Without further restrictions the answer to this question is of
course negative. Indeed, the complex numbers contain all roots of
unity; so there are arbitrarily large finite groups inside $\CC^*$.
Thinking of complex numbers as scalar matrices, we also obtain
arbitrarily large finite groups of $n \times n$-matrices over $\CC$.

The situation changes when certain arithmetic conditions are imposed
on the matrix group $\G$. When all matrices in $\G$ have entries in
the field $\QQ$ rational numbers, Minkowski \cite{hM87} has shown
that the order of $\G$ divides some explicit, and optimal, constant
$M(n)$ depending only on the matrix size $n$. Later, Schur
\cite{iS05} improved on this result by showing that Minkowski's
bound $M(n)$ still works if only the traces of all matrices in $\G$
are required to belong to $\QQ$.

\subsection{} 

The first four sections of this article present full proofs of the
theorems of Schur and Minkowski that depend on very few
prerequisites. These sections follow Schur's approach via character
theory and have been written with a readership of beginning graduate
and advanced undergraduate students in mind. Provided the reader is
willing to accept one simple fact concerning group representations
(Fact 2 in Section~\ref{SS:characters} below), the proofs will be
completely understandable with only a rudimentary knowledge of
linear algebra, group theory (symmetric groups, Sylow's theorem),
and some algebraic number theory (minimal polynomials, Galois groups
of cyclotomic fields). The requisite background material will be
reviewed in Section~\ref{S:tools}.

The material in Section~\ref{S:modp} is new. We show that
Minkowski's original approach used in \cite{hM87} in fact also
yields Schur's theorem \cite{iS05}. Minkowski's method is
conceptually very simple, and it quickly and elegantly explains why
some bound on the order $|\G|$ must exist, even for arbitrary
algebraic number fields, that is, finite extensions of $\QQ$. The
method proceeds by reduction modulo suitably chosen primes and then
using information about the orders of certain classical linear
groups over finite fields. In fact, the general linear group alone
almost suffices; only dealing with the $2$-part of $|\G|$ using this
strategy requires additional information. Since we work over
algebraic number fields, a bit more mathematical background is
assumed in this section.

As of this writing, Schur's theorem first appeared in print exactly
a century ago and Minkowski's goes even further back. In the final
section of this article, we will survey some recent related work of
Collins, Feit and Weisfeiler on finite groups of matrices, in
particular on the so-called Jordan bound. We will also mention two
mysterious coincidences concerning the Minkowski numbers $M(n)$, one
proven but unexplained, the other merely based on experimental
evidence as of now.

\subsection{} 

Minkowski \cite{hM87} proved his remarkable theorem in the course of
his investigation of quadratic forms. Stated in group theoretical
terms, the theorem reads as follows.

\begin{thm}[Minkowski 1887] \label{T:minkowski}
The least common multiple of the orders of all finite groups of $n
\times n$-matrices over $\QQ$ is given by
\begin{equation} \label{E:minkowskibound}
M(n) = \prod_p p^{\left\lfloor \frac{n}{p-1}\right\rfloor +
\left\lfloor \frac{n}{p(p-1)}\right\rfloor + \left\lfloor
\frac{n}{p^2(p-1)}\right\rfloor + \dots}
\end{equation}
\end{thm}

Here, $\left\lfloor \,.\,\right\rfloor$ denotes the greatest integer
less than or equal to $\,.\,$ and $p$ runs over all primes. Note
that if $p > n+1$ then the corresponding factor in the product
equals $1$ and can be omitted. Therefore, \eqref{E:minkowskibound}
is actually a finite product. The first few values of $M(n)$ are:
$$
M(1) = 2^1 = 2\,, \ M(2) = 2^{2+1}\,3^1 = 24\,,\ M(3) = 2^{3+1}\,3^1
= 48\,, \ M(4) = 2^{4+2+1}\,3^{2}\,5^1 = 5760 \ .
$$

\subsection{} \label{SS:p}

For a positive integer $m$ and a prime $p$,  let $m_p$ denote the
$p$-part of $m$, that is, the largest power of $p$ dividing $m$.
Thus, $M(n)_p = p^{\left\lfloor \frac{n}{p-1}\right\rfloor +
\left\lfloor \frac{n}{p(p-1)}\right\rfloor + \left\lfloor
\frac{n}{p^2(p-1)}\right\rfloor + \dots}$. This number can be
written in a more compact form. Indeed, the $p$-part of $m! = 1\cdot
2 \cdot \ldots \cdot m$ is given by
\begin{equation} \label{E:factorial}
(m !)_p = p^{\left\lfloor \frac{m}{p}\right\rfloor + \left\lfloor
\frac{m}{p^2}\right\rfloor + \dots}\ .
\end{equation}
To see this, put $m' = \left\lfloor \frac{m}{p}\right\rfloor$ and
note that $ m! = p\cdot (2p) \cdot \ldots \cdot (m'p) \cdot
(\text{factors not divisible by $p$})$\,. Therefore, $(m !)_p  =
p^{m'}(m' !)_p$ and \eqref{E:factorial} follows by induction. Using
\eqref{E:factorial} we can write
\begin{equation} \label{E:p2}
M(n)_p = p^{\left\lfloor \frac{n}{p-1}\right\rfloor}\left(
\left\lfloor \tfrac{n}{p-1}\right\rfloor !\right)_p\ .
\end{equation}

\subsection{} 
The notation $M(n)$, in the variant $M_n$, was introduced by Schur
in \cite{iS05} to honor Minkowski who had originally denoted the
same number by $\overline{\stackrel{\phantom{
}}{\phantom{.}\!n}}\!\big|$. Relaxing the condition in
Theorem~\ref{T:minkowski} that all matrix entries be rational and
replacing it with the weaker requirement that only the matrix traces
belong to $\QQ$, Schur was able to prove that Minkowski's bound
$M(n)$ still works:

\begin{thm}[Schur 1905] \label{T:schur}
If $\G$ is any finite group of $n \times n$-matrices over $\CC$ such
that $\tr(g) \in \QQ$ holds for all $g \in \G$ then the order of
$\G$ divides $M(n)$.
\end{thm}

Schur's theorem covers a considerably larger class of groups than
Theorem~\ref{T:minkowski}. In \cite{iS05}, the following example of
a group covered by Theorem~\ref{T:schur} but not
Theorem~\ref{T:minkowski} is given.

\begin{example}
Consider the matrices $g = \left( \begin{smallmatrix}  0 & -1 \\ 1 &
0 \end{smallmatrix} \right)$ and $h = \left( \begin{smallmatrix} i &
0 \\ 0 & -i \end{smallmatrix} \right)$, where $i = \sqrt{-1} \in
\CC$. Then $g^2 = h^2 = - 1_{2 \times 2}$ and $gh = \left(
\begin{smallmatrix}  0 & i \\ i & 0 \end{smallmatrix} \right) =
-hg$. Thus $\G = \{ \pm 1_{2 \times 2}, \pm g, \pm h, \pm gh \}$ is
a group of complex $2 \times 2$-matrices of order $8$; it is
isomorphic to the so-called quaternion group $\mathcal{Q}_8$. Note
that the traces of all elements of $\G$ are rational -- they are
either $0$ or $\pm 2$ -- but $\G$ certainly does not consist of
matrices over $\QQ$. In fact, there does not even exist an
invertible complex $2 \times 2$-matrix $a$ such that the matrices $x
= aga^{-1}$ and $y = aha^{-1}$ both have entries in the field $\RR$
of real numbers. To see this, note that $x$ and $y$ both would have
determinant $1$ and trace $0$, as $g$ and $h$ do. A direct
calculation shows that the product matrix $z = xy$ then satisfies
${z_{12}}^2 + {x_{12}}^2 + {y_{12}}^2 = - x_{12}y_{12}\tr(z)$, where
$\,.\,_{12}$ indicates the $(1,2)$-entry of the matrix in question.
However, $\tr(z) = \tr(gh) = 0$. Hence, if $x$ and $y$ are matrices
over $\RR$ then all terms on the left will be zero. But then $1 =
\det(x) = x_{11}x_{22} = -{x_{11}}^2$ which is impossible.
\end{example}

We remark in passing that, for any ``irreducible" finite group $\G$
of complex $n \times n$-matrices,  a  necessary and sufficient
condition for the existence of an invertible complex $n \times
n$-matrix $a$ such that $aga^{-1}$ is real for all $g \in \G$ is
that
$$
\frac{1}{|\G|} \sum_{g \in \G} \tr(g^2) = 1 \ .
$$
The sum on the left is called the \emph{Frobenius-Schur indicator}
of $\G$; see, e.g., Isaacs \cite[Chapter 4]{mI76}. The group $\G =
\mathcal{Q}_8$ in the example above has Frobenius-Schur indicator
$-1$.

\subsection{} 

The proof of Theorems~\ref{T:minkowski} and  \ref{T:schur} to be
given in Section~\ref{S:proof} below proceeds by first exhibiting
sufficiently large groups of rational (in fact, integer) matrices
showing that the least common multiple of the orders of all finite
groups of $n \times n$-matrices over $\QQ$ must be at least equal to
$M(n)$. Thereafter, we may concentrate on Theorem~\ref{T:schur}
which in particular implies that the least common multiple in
Theorem~\ref{T:minkowski} does not exceed $M(n)$. Apart from
updating terminology and notation to current usage and adding more
generous details to the exposition, we have followed Schur's
original approach in \cite{iS05} quite closely. For a proof of
Schur's theorem using slightly more sophisticated tools from
representation theory, see Isaacs \cite[Theorem 14.19]{mI76}.
Stronger results are presented in Feit \cite{wF97}.

\subsection{} 

This article is dedicated to our friend and colleague Don Passman.
Don's contributions to group theory and ring theory in general and
his expository masterpieces \cite{dP71}, \cite{dP77} in particular
have profoundly influenced our own work. In the course of various
collaborations with Don, we have both benefitted from his deep
insights and his generosity in sharing ideas.

\begin{notat}
Throughout, $\GL_n(R)$ will denote the group of all invertible $n
\times n$-matrices over the commutative ring $R$. Recall that a
matrix over $R$ is invertible if and only if its determinant is an
invertible element of $R$.
\end{notat}


\section{Large groups of integer matrices} \label{S:large}

The principal goal of this section is to construct certain groups of
$n \times n$-matrices over $\ZZ$ such that the least common multiple
of their orders equals the Minkowski bound $M(n)$ in
\eqref{E:minkowskibound}. This will then allow us to give a
reformulation of the core of Theorem~\ref{T:schur}.

\subsection{Construction of groups} \label{SS:construction}
The main building blocks of the construction will be the symmetric
groups $\Sy_r$ for various $r$. Recall that $\Sy_r$ consists of all
permutations of $\{1,\dots,r\}$ and has order $r!$\,.

\begin{prop} \label{P:minkowski}
Let  $a$, $m$ and $n$ be positive integers with $am \le n$.
Then $\GL_n(\ZZ)$ has a subgroup $\G$ of order $|\G| =
(m+1)!^{\,a}\, a!$\ .
\end{prop}

\begin{proof}
If we can realize $\G$ inside $\GL_{am}(\ZZ)$ then we can view $\G$
as a subgroup of $\GL_n(\ZZ)$ via
$$
\G \into \GL_{am}(\ZZ) \cong
\left(
      \begin{tabular}{cc}
        \fbox{\rule[-3mm]{0mm}{11mm}$\GL_{am}(\ZZ)$} &  \\
       & $\begin{matrix} 1 & & \\ & \ddots & \\ & & 1 \end{matrix}$
      \end{tabular}
\right)
\into \GL_{n}(\ZZ) \ .
$$
Therefore, we may assume that $n = am$. Think of the rows of any $n
\times n$-matrix as partitioned into $a$ blocks of $m$ adjacent
rows, and similarly for the columns. Now consider all matrices in
$\GL_n(\ZZ)$ that have exactly one $m \times m$-identity matrix
$1_{m \times m}$ in each block of rows and each block of columns and
$0$s elsewhere; these are special permutation matrices. In fact, the
collection of all these matrices forms a subgroup  $\Pi \subseteq
\GL_n(\ZZ)$ that is isomorphic to the symmetric group $\Sy_a$\,:
$$
\Sy_a \cong \Pi = \left\{ \left(
      \begin{tabular}{c|c|c|c|c|c}
        & $1_{m \times m}$ & $\dots$ & & & \\
        \cline{1-6}
        & & $\dots$ & & $1_{m \times m}$ & \\
        \cline{1-6}
        & & $\ddots$ & & & \\
        \cline{1-6}
        & & $\dots$& $1_{m \times m}$ & &  \\
      \end{tabular}
\right)\right\}
\into \GL_{n}(\ZZ) \ .
$$
Next, we turn to the symmetric group $\Sy_{m+1}$. This group acts on the lattice $\ZZ^{m+1}$
by permuting
its canonical basis $e_1 = (1,0,\dots,0), \dots , e_{m+1} =
(0,\dots,0,1)$ via $\sigma(e_i) = e_{\sigma(i)}$. Note that this
action maps the following sublattice to itself:
$$
A_m = \{ (z_1,\dots,z_{m+1}) \in \ZZ^{m+1} \mid \sum_i z_i = 0 \} \cong \ZZ^m
$$
(The notation $A_m$ comes from the theory of root systems;
cf.~\cite{nB68}.) Thus, fixing some $\ZZ$-basis of $A_m$, each
permutation $\sigma \in \Sy_{m+1}$ yields a matrix
$\widetilde{\sigma} \in \GL_m(\ZZ)$. It is easy to see that the map
$\sigma \mapsto \widetilde{\sigma}$ is an injective group
homomorphism $\Sy_{m+1} \to \GL_m(\ZZ)$. Stringing each $a$-tuple
$(\widetilde{\sigma_1},\dots,\widetilde{\sigma_a})$ along the
diagonal in $\GL_n(\ZZ)$ we obtain a subgroup $\Delta \subseteq
\GL_n(\ZZ)$ that is isomorphic to $\Sy_{m+1}^a$\,:
$$
\Sy_{m+1}^a =
\underbrace{ \Sy_{m+1} \times \dots \times \Sy_{m+1}
}_{\mbox{$a$ factors}} \cong \Delta = \left\{ \left(
      \begin{tabular}{cccc}
        \fbox{\rule[-2mm]{0mm}{6mm}$\widetilde{\sigma_1}$} & & & \\
        & \fbox{\rule[-2mm]{0mm}{6mm}$\widetilde{\sigma_2}$} & & \\
        & & $\ddots$ &  \\
        & & & \fbox{\rule[-2mm]{0mm}{6mm}$\widetilde{\sigma_a}$}  \\
      \end{tabular}
\right)\right\}
\into \GL_{n}(\ZZ) \ .
$$
The subgroup $\Pi$ of $\GL_n(\ZZ)$ constructed earlier has only the
identity matrix in common with $\Delta$. Moreover, conjugating a
matrix in $\Delta$ with a matrix from $\Pi$ simply permutes the
$\widetilde{\sigma_i}$-blocks along the diagonal. Therefore,
defining $\G$ to be the subgroup of $\GL_n(\ZZ)$ that is generated
by $\Pi$ and $\Delta$, we obtain
$$|\G| = |\Delta|\,|\Pi| =
(m+1)!^{\,a}\, a!\ ,
$$
as desired.
\end{proof}


Now fix a prime $p \le n+1$. Taking $m = p-1$ and $a = \left\lfloor
\frac{n}{p-1}\right\rfloor$ in Proposition~\ref{P:minkowski} we
obtain a subgroup $\G$ of $\GL_n(\ZZ)$ of order $p!^{\,a}\, a !$; so
$|\G|_p = p^a(a!)_p$. In view of \eqref{E:p2}, this says that
$|\G|_p = M(n)_p$. Letting $p$ range over all primes $\le n+1$, we
have exhibited a collection of subgroups of $\GL_n(\ZZ)$ such that
the least common multiple of their orders is $M(n)$\,.

\subsection{Reformulation of Theorem~\ref{T:schur}} \label{SS:reduction}

Let $\G \subseteq \GL_n(\CC)$ be as in Theorem~\ref{T:schur}. Our
goal is to show that, for all primes $p$, the $p$-part $|\G|_p$
divides $M(n)_p = p^{a}\, (a !)_p$ with $a = \left\lfloor
\frac{n}{p-1}\right\rfloor$ as in \eqref{E:p2}. Now Sylow's Theorem
tells us that $\G$ has subgroups of order $|\G|_p$, the so-called
Sylow $p$-subgroups of $\G$. Replacing $\G$ by one of its Sylow
$p$-subgroups, the issue becomes to show that $|\G|$ divides
$p^{a}a!$\,. Therefore, in order to prove Theorem~\ref{T:schur}, and
thereby complete the proof of Theorem~\ref{T:minkowski}, it suffices
to establish the following proposition.

\begin{prop} \label{P:reduction}
Let $\G$ be finite subgroup of $\GL_n(\CC)$ whose order is a
$p$-power for some prime $p$ and such that $\ \tr(g) \in \QQ$ holds
for all $g \in \G$. Then $|\G|$ divides $p^{a}a!$ with $a =
\left\lfloor \frac{n}{p-1}\right\rfloor$.
\end{prop}


\section{Tools for the proof} \label{S:tools}

The proof of Proposition~\ref{P:reduction} will depend on three
ingredients: a lemma to narrow down the possible trace values, some
basic facts on characters of group representations, and an
observation concerning the familiar Vandermonde matrix. We will
discuss each of these topics in turn.

\subsection{Traces} \label{SS:traces}
This section uses a small amount of algebraic number theory. The
book \cite{gJ96} by Janusz is a good background reference.

Besides the usual matrix traces, we will use a notion of trace that
is associated with field extensions. Specifically, let $K/F$ be a
finite Galois extension with Galois group $\Gamma = \Gal(K/F)$. Then
the trace $\Tr_{K/F} \colon K \to F$ is defined by
$\Tr_{K/F}(\alpha) = \sum_{\gamma \in \Gamma} \gamma(\alpha)$ for
$\alpha \in K$. If $x^{m} + cx^{m-1}+\dots$ is the minimal
polynomial of $\alpha$ over $F$ then
\begin{equation} \label{E:galoistr}
\Tr_{K/F}(\alpha) = - \frac{|\Gamma|}{m} \cdot c \ .
\end{equation}
This follows from the fact that the minimal polynomial of $\alpha$
is equal to $\prod_{i=1}^m(x - \alpha_i)$, where $\{\alpha_i\}_1^m$
are the distinct Galois conjugates $\gamma(\alpha)$ with $\gamma \in
\Gamma$. We will only be concerned with the special case where $F =
\QQ$ and $K = \QQ(e^{2\pi i/p^r})$ with $p$ prime. The Galois group
of $\QQ(e^{2\pi i/p^r})/\QQ$ is isomorphic to the group of units
$\left(\ZZ/p^r\ZZ\right)^*$ of the ring $\ZZ/p^r\ZZ$; its order is
$\varphi(p^r) = p^{r-1}(p-1)$.

\begin{lem} \label{L:algebraic}
Let $g \in \GL_n(\CC)$ be a matrix whose order is a power of $p$ and
such that $\tr(g) \in \QQ$. Then $\tr(g)$ must be one of the values
$\{n, n-p, n-2p, \dots , n-ap \}$, where $a = \left\lfloor
\frac{n}{p-1}\right\rfloor$. Moreover, $\tr(g) = n$ holds only for
$g=1_{n \times n}$.
\end{lem}

\begin{proof}
By hypothesis, $g^{p^r} = 1_{n \times n}$ for some $r$. Let
$\e_1,\dots,\e_n$ denote the eigenvalues of $g$; they are all powers
of $\z = e^{2\pi i/p^r}$. Hence, $\tr(g) = \sum_i \e_i$ belongs to
the subring $\ZZ[\z] \subseteq \CC$ while also being rational, by
hypothesis. This implies that  $\tr(g)$ is actually an integer; see
\cite[Section I.2]{gJ96}. Furthermore, by the triangle inequality,
$|\tr(g)| \le \sum_i |\e_i| = n$ and $\le$ is equality if and only
if all $\e_i$ are the same, that is, $g$ is a scalar matrix. In
particular, $\tr(g) = n$ holds only for $g=1_{n \times n}$.

Let $\p = (\z - 1)$ denote the ideal of $\ZZ[\z]$ that is generated
by the element $\z - 1$. So $\z \equiv 1 \mod \p$, and hence all
$\e_i \equiv 1 \mod \p$ and $\tr(g) \equiv n \mod \p$. Therefore,
$\tr(g) - n \in \p \cap \ZZ = (p)$; see \cite[Theorem I.10.1]{gJ96}
for the last equality. Since we have already shown that $\tr(g) \le
n$,  we conclude that $\tr(g) = n - pt$ for some non-negative
integer $t$. It remains to show that $t \le \frac{n}{p-1}$ or,
equivalently,
$$
\tr(g) \ge -\frac{n}{p-1} \ .
$$
To this end, consider the Galois extension $\QQ(\z)/\QQ$ and its
trace $\Tr_{\QQ(\z)/\QQ}$. The minimal polynomial over $\QQ$ of a
root of unity of order $p^s >1$ is given by $x^{p^{s-1}(p-1)} +
x^{p^{s-1}(p-2)} + \ldots + 1$ (\cite[Theorem I.10.1]{gJ96} again).
Therefore, equation \eqref{E:galoistr} yields
$$
\Tr_{\QQ(\z)/\QQ}(\e_i) = \begin{cases} \varphi(p^r) & \text{if $\e_i = 1$\,,} \\
- p^{r-1} & \text{if $\e_i$ has order $p$\,,} \\
0 & \text{otherwise.}
\end{cases}
$$
Put $n_0 = \#\{i \mid \e_i = 1 \}$ and $n_1 = \#\{i \mid
\text{$\e_i$ has order $p$} \}$; so $0 \le n_i \le n$. Using the
fact that $\tr(g) \in \QQ$ we obtain
$$
\varphi(p^r) \tr(g) = \Tr_{\QQ(\z)/\QQ}(\tr(g)) =
\sum_i \Tr_{\QQ(\z)/\QQ}(\e_i) = \varphi(p^r) n_0 - p^{r-1} n_1 \ .
$$
Hence, $\tr(g) = n_0 - \frac{n_1}{p-1} \ge -\frac{n}{p-1}$, as
desired.
\end{proof}

\subsection{Characters} \label{SS:characters}
A \emph{complex representation} of a group $\G$ is a homomorphism
$\rho \colon \G \to \GL(V)$ for some $\CC$-vector space $V$. If $n =
\dim_{\CC}V$ then we may identify $\GL(V)$ with $\GL_n(\CC)$; the
integer $n$ is called the \emph{degree} of the representation
$\rho$. The \emph{character} $\chi = \chi_{\rho}$ of $\rho$ is the
complex-valued function on $\G$ that is given by $\chi(g) =
\tr(\rho(g))$ for $g \in \G$. \medskip
\begin{quote}
\textbf{Fact 1} \emph{The sum $\sum_{g \in \G} \chi(g)$ is always an
integer that is divisible by  $|\G|$.}
\end{quote}\medskip
To see this, consider the linear operator $e_{\rho} \in
\End_{\CC}(V) \cong \M_n(\CC)$ that is defined by $e_{\rho} =
\frac{1}{|\G|} \sum_{g \in \G} \rho(g)$. Note that $\rho(g) e_{\rho}
= e_{\rho}$ holds for all $g \in \G$, because multiplication with
$\rho(g)$ simply permutes the summands of $e_{\rho}$. Hence,
$e_{\rho}$ is an idempotent operator: $e_{\rho}^2 = e_\rho$.
Therefore, the trace of $e_\rho$ is equal to the rank of $e_\rho$:
$\tr(e_\rho) = \dim_{\CC}e_{\rho}(V)$. On the other hand,
$\tr(e_\rho) = \frac{1}{|\G|} \sum_{g \in \G} \tr(\rho(g)) =
\frac{1}{|\G|} \sum_{g \in \G} \chi(g)$. This proves Fact 1. We
remark that Fact 1 is a special case of the so-called
\emph{orthogonality relations} of characters.\medskip
\begin{quote}
\textbf{Fact 2} \emph{The product of any two characters of $\G$ is
again a character of $\G$. In particular, all powers $\chi^s \ (s
\ge 0)$ of a character $\chi$ are also characters of $\G$.}
\end{quote}\medskip
Here, the $0^\text{th}$ power $\chi^0$ is the constant function with
value $1$; it is the character of the so-called trivial
representation $\G \to \CC^* = \GL_1(\CC)$ sending every $g \in \G$
to $1$. In order to show that the product of two characters,
$\chi_{\rho}$ and $\chi_{\rho'}$, is itself a character, one needs
to construct a complex representation of $\G$ whose character is
$\chi_{\rho}\cdot\chi_{\rho'}$. This is achieved by the so-called
tensor product $\rho \otimes \rho'$ of the representations $\rho$
and $\rho'$, a complex representation of degree equal to
$\deg\rho\cdot\deg\rho'$ for whose detailed construction the reader
is referred to Isaacs \cite[Chapter 4]{mI76} or any other text on
group representation theory. More generally, tensor products of
representations can be defined for \emph{Hopf algebras}; they form
an important aspect of the current investigation of \emph{quantum
groups}.

\subsection{Vandermonde matrix} \label{SS:vandermonde}
Given a collection $z_0,\dots,z_a$ of elements in some commutative ring $R$
(later we will take $R = \ZZ$), form the
familiar Vandermonde  matrix
$$
V  = \begin{pmatrix}
1 & z_0 & z_0^2 & \dots & z_0^{a} \\
1 & z_1 & z_1^2 & \dots & z_1^{a} \\
\vdots & \vdots & \vdots & \ddots & \vdots \\
1 & z_a & z_a^2 & \dots & z_a^{a}
\end{pmatrix}\ .
$$
We will exhibit a matrix $E$ over $R$ so that the
matrix product $V\cdot E$ is diagonal:
\begin{equation} \label{E:vandermonde}
V \cdot E = \diag\left(\prod_{\substack{0 \le s \le a \\ s \neq 0}} z_0 - z_s,
\prod_{\substack{0 \le s \le a \\ s \neq 1}} z_1 - z_s, \dots ,
\prod_{\substack{0 \le s \le a \\ s \neq a}} z_a - z_s \right) \ .
\end{equation}
To this end, let $e_s = e_s(x_1,\dots,x_a)$ denote the $s^\text{th}$
elementary symmetric function in the commuting variables
$x_1,\dots,x_a$\,. These functions can be defined by
\begin{equation} \label{E:symmetric}
\prod_{i=1}^a (x - x_i) = \sum_{s=0}^a x^{s}(-1)^{a-s} e_{a-s}\ ,
\end{equation}
where $x$ is an additional commuting variable. Explicitly, $e_s =
\sum_I \prod_{i \in I} x_i$, where $I$ runs over all subsets $I
\subseteq \{1,\dots,a\}$ with $|I| = s$. Specializing $x$ to
$z_{t'}$ and $(x_1,\dots,x_a)$ to
$(z_0,\dots,\widehat{z_t},\dots,z_a)$, where $\widehat{z_t}$ signals
that $z_t$ has been deleted from the list, and defining
$$
E = \left((-1)^{a-s} e_{a-s}(z_0,\dots,\widehat{z_t},\dots,z_a)\right)_{s,t=0,\dots,a}
$$
equation \eqref{E:symmetric} becomes the desired equation \eqref{E:vandermonde}.


\section{Schur's proof of Theorems~\ref{T:minkowski} and~\ref{T:schur}} \label{S:proof}

It remains to prove Proposition~\ref{P:reduction}. So fix a prime
$p$ and let $\G$ be finite subgroup of $\GL_n(\CC)$ whose order
$|\G|$ is a power of $p$. We assume that $\tr(g) \in \QQ$ holds for
all $g \in \G$. Since the order of each $g$ divides $|\G|$,
Lemma~\ref{L:algebraic} implies that the traces $\tr(g)$ can only
take the  values
$$
z_t = n - pt \quad\text{with $0 \le t \le a = \left\lfloor
\tfrac{n}{p-1}\right\rfloor$.}
$$
Put $m_t = \#\{ g \in \G \mid \tr(g) = z_t \}$; so  $m_0 = 1$ by
Lemma~\ref{L:algebraic}. Proposition~\ref{P:reduction} is the case
$t=0$ of the following

\begin{claim}
For all $0 \le t \le a$, the order $|\G|$ divides the product
$\displaystyle m_t p^a \prod_{\substack{0 \le s \le a \\ s \neq t}}
s-t$.
\end{claim}

\noindent  To prove this, note that the inclusion $\G \into
\GL_n(\CC)$ is a complex representation of $\G$ with character
$\chi(g) = \tr(g)$. Therefore, it follows from Facts 1 and 2 above
that, for each non-negative integer $s$, the sum $\sum_{g \in \G}
\tr(g)^s$ is an integer that is divisible by $|\G|$. In other words,
$\sum_{t=0}^a m_t z_t^s \equiv 0 \mod |\G|$ or, in matrix form,
\begin{equation} \label{E:first}
(m_0,\dots,m_a) \cdot V \equiv (0,\dots,0) \mod |\G|\ ,
\end{equation}
where $V = \left(z_t^s\right)_{t,s=0,\dots,a}$ is the Vandermonde
matrix, as in \S\ref{SS:vandermonde}. Multiplying both sides of
equation \eqref{E:first} with the matrix $E$ constructed in
\S\ref{SS:vandermonde}, we deduce from equation
\eqref{E:vandermonde}  that
$$
m_t \prod_{\substack{0 \le s \le a \\ s \neq t}} z_t - z_s \equiv 0
\mod |\G|
$$
holds for all $0 \le t \le a$. Since $z_t - z_s = p(s-t)$, this is exactly what the claim states.
This completes the proof of Proposition~\ref{P:reduction}, and hence
Theorems~\ref{T:minkowski} and~\ref{T:schur} are proved as well.


\section{Minkowski's reduction method} \label{S:modp}

Minkowski's original proof of Theorem~\ref{T:minkowski} is quite
different from Schur's. The essential tool are reduction
homomorphisms to the general linear group over certain finite
fields. The reduction method applies to algebraic number fields $K$,
that is, finite extensions of $\QQ$, and very quickly yields rough
bounds for the orders of all finite subgroups $\G \subseteq
\GL_n(K)$; see Proposition~\ref{P:p'} below. In fact, subgroups $\G
\subseteq \GL_n(\CC)$ satisfying only $\tr(g) \in K$ for all $g \in
\G$ can also be treated by this strategy due to the fact that linear
groups over finite fields can be realized over the subfield
generated by the traces; see Lemma~\ref{L:galois}. A sharp bound for
the $2'$-part of $|\G|$ can be easily deduced in this way from the
well-known order of the general linear group over a finite field
together with some elementary number theoretic observations; see
Proposition~\ref{P:oddell}. The $2$-part of $|\G|$ requires
additional information concerning certain classical groups
associated to hermitian or skew-hermitian forms. This will be
explained in \S\S~\ref{SS:bilinear} and \ref{SS:congruence3} below.

As usual, the field with $q$ elements will be denoted by $\FFq$. We
will also occasionally write the $p$-part of an integer $m$ as $m_p
= p^{v_p(m)}$, and $m_{p'}$ will denote the $p'$-part of $m$; so
$m_{p'} = m/m_p$.

\subsection{The general linear group over finite fields} \label{SS:linear}

It is well-known and easy to see that $\GL_n(\FFq)$ has order
$\prod_{i=0}^{n-1} (q^n - q^i)$; cf., e.g., Rotman \cite[Theorem
8.5]{jjRot95}. Thus, if $q = p^f$ then
\begin{equation} \label{E:GLp'}
|\GL_n(\FFq)|_{p'} = \prod_{i=1}^n (q^i - 1) \ .
\end{equation}

\begin{lem} \label{L:Q}
Let $\ell$ be an odd prime. There are infinitely many primes $p$
such that
$$
|\GL_n(\FF_{\!p^f})|_{\ell} = \ell^{\left( 1 +
v_{\ell}(f)\right)\left\lfloor n/\tau \right\rfloor}
\left(\left\lfloor n/\tau \right\rfloor !\right)_{\ell}
$$
holds for all positive integers $n$ and $f$, where $\tau =
\frac{\ell-1}{(\ell - 1,f)}$.
\end{lem}

\begin{proof}
We use the fact that, for odd primes $\ell$, the group of units
$(\ZZ/\ell^s\ZZ)^*$ of the ring $\ZZ/\ell^s\ZZ$ is cyclic of order
$\varphi(\ell^s) = \ell^{s-1}(\ell-1)$. Any integer whose residue
class modulo $\ell^2$ generates $(\ZZ/\ell^2\ZZ)^*$ will also
generate the units modulo all powers $\ell^s$; see \cite[proof of
Theorem 2 on p.~43]{kImR93}. Moreover, by Dirichlet's theorem on
primes in arithmetic progression (e.g., \cite[p.~61]{jpS73}), the
residue class modulo $\ell^2$ of any generator of
$(\ZZ/\ell^2\ZZ)^*$ contains infinitely many primes $p$. Let $p$ be
one of these primes. Then $p$ has order $\varphi(\ell^s)$ in
$(\ZZ/\ell^s\ZZ)^*$; so $p^i \equiv 1 \mod \ell^s$ if and only if
$i$ is divisible by $\varphi(\ell^s)$. In other words, $\ell$
divides $p^i - 1$ if and only if $\ell-1$ divides $i$ and, in this
case,
\begin{equation*}
(p^i - 1)_{\ell} =  \ell \left( \tfrac{i}{\ell-1} \right)_{\ell} \ .
\end{equation*}
Now put $q = p^f$. Then $\ell$ divides $q^i - 1$ if and only if
$\tau$ divides $i$ and, in this case, $(q^i - 1)_{\ell} = \ell\,
f_\ell\, (i/\tau)_\ell $. For $1\le i \le n$, this applies to $i =
\tau,2\tau,\dots,\alpha\tau$, where $\alpha = \left\lfloor n/\tau
\right\rfloor$. Thus, $|\GL_n(\FF_{\!p^f})|_{\ell} = \prod_{i=1}^n
(q^i-1)_{\ell} = \left( \ell\, f_\ell \right)^\alpha
(\alpha!)_{\ell}$\,, which proves the lemma.
\end{proof}

We remark that, for $f=1$, the expression $\ell^{\left( 1 +
v_p(f)\right)\left\lfloor n/\tau \right\rfloor} \left(\left\lfloor
n/\tau \right\rfloor !\right)_{\ell}$ in  Lemma~\ref{L:Q} is
identical with the $\ell$-part of the Minkowski bound $M(n)$; see
equation \eqref{E:p2}. Thus, for an odd prime $\ell$,
\begin{equation} \label{E:p1}
|\GL_n(\FF_{\!p})|_{\ell} = M(n)_\ell
\end{equation}
holds for infinitely many primes $p$. Lemma~\ref{L:Q} fails for the
prime $\ell = 2$, because the linear group is too big. For example,
for all odd primes $p$, $|\GL_2(\FF_{\!p})|_{2} = (p-1)_2(p^2-1)_2$
is divisible by $16$ while $M(2)_2 = 8$.

\begin{lem} \label{L:galois}
Let $\G$ be a finite subgroup of $\GL_n(\FFq)$, where $q = p^f$.
Assume that $p$ does not divide $|\G|$ and that $p > n$. If all $g
\in \G$ satisfy $\tr(g) \in F$ for some subfield $F \subseteq \FFq$
then $\G$ is conjugate to a subgroup of $\GL_n(F)$.
\end{lem}

\begin{proof}
Let $\k = F^{\text{\rm alg}}$ denote an algebraic closure of $F$
with $\FFq \subseteq \k$, and let  $\sigma$ denote the canonical
topological generator of $\Gal(\k/F) \cong \widehat{\ZZ}$. Then
$\sigma$ acts on $\GL_n(\k)$ by $\left(g_{i,j}\right)^{\sigma}_{n
\times n} = \left(g_{i,j}^{\sigma}\right)_{n \times n}$. By our
hypothesis on traces, the map $\G \to \GL_n(\k)$, $g \mapsto
g^{\sigma}$, is a $\k$-representations of $\G$ having the same
character as the inclusion $\G \hookrightarrow \GL_n(\k)$. Since
both representations are semisimple, by Maschke's theorem, they are
isomorphic. (The proof of \cite[\S~12, Proposition 3]{nB73} works in
characteristic $p > n$.) Thus, there exists a matrix $u \in
\GL_n(\k)$ such that $u g u^{-1} = g^\sigma$ holds for all $g \in
\G$. By Lang's theorem \cite{sL56}, we can write $u =
v^{\sigma}v^{-1}$ for some $v \in \GL_n(\k)$. Thus, each $v^{-1}gv$
is fixed by $\sigma$, and hence it belongs to $\GL_n(F)$. By the
Noether-Deuring Theorem (e.g., Curtis-Reiner \cite[p.~139]{cCiR81}),
we may replace $v$ by a matrix in $\GL_n(\FFq)$, proving the lemma.
\end{proof}

\begin{remarks}
(a) Lang's theorem is a much more general result than what is
actually needed for the proof of Lemma~\ref{L:galois}; see, e.g.,
Borel \cite[Corollary 16.4]{aB91}. Indeed, we only invoke the
theorem for the algebraic group $\GL_n$ and, in this case, it is a
special case of Speiser's version of Hilbert's Theorem 90: the
Galois cohomology set $H^1(F,\GL_n)$ is trivial for every field $F$;
cf. Serre \cite[Proposition X.3]{jpS79} or Knus et.~al.~\cite[Remark
29.3]{KMRT}. For a finite field $F$, triviality of $H^1(F,\GL_n)$
amounts to the desired fact that every $u \in \GL_n(F^{\text{alg}})$
can be written as $u = v^{\sigma}v^{-1}$, where $\sigma$ is the
Frobenius generator of $\Gal(F^{\text{alg}}/F)$; see \cite[Exercise
2 on p.~442]{KMRT}.

(b) It follows from (a) that $H^1(\FFq,\PGL_n)$ is trivial as well:
every $U \in \PGL_n(\FFq^{\text{alg}}) = \GL_n(\FFq^{\text{alg}})/(
\FFq^{\text{alg}} )^*$ can be written as $U = V^{\sigma}V^{-1}$ for
some $V \in \PGL_n(\FFq^{\text{alg}})$. Moreover, triviality of
$H^1(\FFq,\PGL_n)$ is equivalent to Wedderburn's commutativity
theorem for finite division rings; see \cite[Proposition X.8]{jpS79}
or \cite[p.~396]{KMRT}. For an alternative proof of a version of
Lemma~\ref{L:galois} based on Wedderburn's commutativity theorem,
see Isaacs \cite[Theorem 9.14]{mI76}. Incidentally, Wedderburn's
article \cite{jW05} appeared in 1905, as did Schur's, and Speiser's
generalization of Hilbert's Theorem 90 appeared in 1919 \cite[Satz
1]{aS19}. None of this was available to Minkowski when \cite{hM87}
was written.
\end{remarks}

\subsection{The reduction map} \label{SS:congruence}

Throughout this section, $K$ will denote an algebraic number field
and $\G$ will be a finite subgroup of $\GL_n(K)$. Furthermore, $\cO
= \cO_K$ will denote the ring of algebraic integers in $K$.

Put $L = \sum_{g \in \G} g\cdot \cO^n \subset K^n$; this is a
$\G$-stable finitely generated $\cO$-submodule of $K^n$. If $\cO$ is
a principal ideal domain (or, put differently, $K$ has class number
$1$) then the theory of modules over PIDs tells us that $L$ is
isomorphic to $\cO^n$; see, e.g., Jacobson \cite[Section 3.8]{nJ85}.
Therefore:
\begin{equation*} \label{E:pid}
\text{\emph{If $\cO = \cO_K$ is a PID then $\G$ is conjugate in
$\GL_n(K)$ to a subgroup of $\GL_n(\cO)$.}}
\end{equation*}
For $K = \QQ$, for example, this says that every finite subgroup of
$\GL_n(\QQ)$ can be conjugated into $\GL_n(\ZZ)$. This explains why
it was enough to look at integer matrices rather than matrices over
$\QQ$ in Section \ref{S:large}.

In general, $\cO$ is a Dedekind domain and the foregoing applies
``locally": for every prime ideal $\p$ of $\cO$, the localization
$\cO_{\p}$ is a PID; see Jacobson \cite[Section 10.2]{nJ89}.
Consequently, as above, we may conclude that $\G$ is conjugate in
$\GL_n(K)$ to a subgroup of $\GL_n(\cO_{\p})$, and hence we may
assume that $\G \subseteq \GL_n(\cO_{\p})$ after replacing $\G$ by a
conjugate. In fact, except for finitely many primes of $\cO$, the
group $\G$ is actually contained in $\GL_n(\cO_{\p})$ at the outset:
if $a \in \cO$ is a common denominator for all matrix entries of all
elements of the original $\G$ then $\G \subseteq \GL_n(\cO[1/a])$;
so any prime $\p$ not containing $a$ will do. Now let $\p \neq 0$
and put $(p) = \p \cap \ZZ$. Then $\cO/\p$ is a finite field of
characteristic $p$. The number of elements of $\cO/\p$ is often
called the \emph{absolute} or \emph{counting norm} of $\p$; it will
be denoted by $\N(\p)$. Thus,
$$
\cO/\p \cong \FF_{\!\N(\p)} \quad\text{and}\quad \N(\p) = p^f \ ,
$$
where $f = f(\p/\QQ)$ is the relative degree of $\p$ over $\QQ$.
Reduction of all matrix entries modulo the maximal ideal $\p
\cO_{\p}$ of $\cO_{\p}$ gives a homomorphism
\begin{equation} \label{E:reduction}
\GL_n(\cO_{\p}) \to \GL_n(\FF_{\!\N(\p)})\ ,
\end{equation}
because $\cO_{\p}/\p \cO_{\p} \cong \cO/\p$. The following lemma is
well-known. Only the first assertion will be needed later; the
second is included for its own sake. Recall that, since $\cO_{\p}$
is a local PID, its non-zero ideals are exactly the powers of the
maximal ideal $\p \cO_{\p}$. The \emph{ramification index} of $\p$
over $\QQ$ is the power $e$ such that $p\cO_{\p} = \p^e \cO_{\p}$.

\begin{lem} \label{L:p'}
The kernel of the reduction homomorphism \eqref{E:reduction} has at
most $p$-torsion. In fact, any torsion element $g$ in the kernel
satisfies $g^{p^i} = 1_{n\times n}$ for some $p^i \le ep/(p-1)$.
\end{lem}

\begin{proof}
For each $g \in \GL_n(\cO_{\p})$, define $d(g) = \sup\{ m \mid g -
1_{n\times n} \in \M_n(\p^m \cO_{\p}) \}$; so $d(g) = \infty$ if and
only if $g = 1_{n\times n}$ and $d(g) > 0$ if and only if $g$
belongs to the kernel of \eqref{E:reduction}. Now assume that $0 < d
= d(g) < \infty$ and
write $g = 1_{n\times n} + \pi^dh$, where $\pi$ is a generator of
the ideal $\p \cO_{\p}$ and $h \in \M_n(\cO_{\p}) \setminus \M_n(\p
\cO_{\p})$. Then $g^{r} = 1_{n\times n} + \pi^{d}(r h + s)$ with $s
= \sum_{i=2}^{\ell} \binom{r}{i}\pi^{d(i-1)}h^i \in \M_n(\p
\cO_{\p})$. If $(r,p) = 1$ then $r h + s \notin \M_n(\p \cO_{\p})$
and so $g^{r} \neq 1_{n\times n}$. This shows that the kernel of
\eqref{E:reduction} has at most $p$-torsion.

We claim that any $g \in \GL_n(\cO_{\p})$ with $d = d(g) >0$
satisfies $d(g^p) \ge \min\{e+d,pd \}$, and $d(g^p) = e+d$ if $pd >
e+d$. Indeed, we may assume that $d < \infty$. Writing $g =
1_{n\times n} + \pi^dh$ be as above, we obtain $g^{p} = 1_{n\times
n} + \pi^{dp}h^p + t$ with $t = \sum_{i=1}^{p-1}
\binom{p}{i}\pi^{di}h^i$. Since $p$ divides all binomial
coefficients $\binom{p}{i}$ occurring in $t$, we have $t \in
\M_n(\p^{e+d} \cO_{\p}) \setminus \M_n(\p^{e+d+1} \cO_{\p})$. The
claim follows from this. We conclude in particular that $g^p \neq
1_{n\times n}$ if $\infty > (p-1)d > e$.

Now assume that $g \in \GL_n(\cO_{\p})$ is a torsion-element with $0
< d(g) < \infty$. Then $g^{p^i} = 1_{n\times n}$ for some positive
integer $i$. If $i$ is chosen minimal then our observations in the
previous paragraph imply that $e \ge (p-1)d(g^{p^{i-1}}) \ge
(p-1)p^{i-1}d(g)$. Hence, $p^i \le ep/(p-1)$ which proves our second
assertion.
\end{proof}

The above proof also shows that if $m(p-1) > e$ then there is no
non-trivial torsion in the kernel of the homomorphism
$\GL_n(\cO_{\p}) \to \GL_n(\cO/\p^m)$ that is defined by reduction
of all matrix entries modulo $\p^m \cO_{\p}$.

\begin{example} \label{EX:Q1}
Let $K = \QQ$. Then $\p = (p)$ and $e=1$. Thus, in Lemma~\ref{L:p'},
we must have $i=0$
 when $p$ is an odd prime, and $i \le 1$ when
$p=2$. In other words, the kernel of the reduction map
$\GL_n(\ZZ_{(p)}) \to \GL_n(\FF_{\!p})$ is torsion-free for odd $p$.
For $p=2$, the only non-trivial torsion possible is order $2$. The
kernel of $\GL_n(\ZZ_{(2)}) \to \GL_n(\ZZ/4\ZZ)$ is torsion-free.
\end{example}

The first assertion of Lemma~\ref{L:p'} implies that the $p'$-part
$|\G|_{p'}$ of the order of $\G$ divides
$|\GL_n(\FF_{\N(\p)})|_{p'}$. In view of equation \eqref{E:GLp'},
this yields  the following proposition.

\begin{prop} \label{P:p'}
Let $\G$ be a finite subgroup of $\GL_n(K)$, where $K$ is an
algebraic number field. Then, for each non-zero prime $\p$ of
$\cO_K$ lying over $p \in \ZZ$, $|\G|_{p'}$ divides $\prod_{i=1}^n
(\N(\p)^i - 1)$\,.
\end{prop}

Applying Proposition~\ref{P:p'} with any two choices of $\p$ lying
over different rational primes yields a bound for the order of $\G$.
Moreover, Proposition~\ref{P:p'} comes close to establishing the
Minkowski bound $M(n)$ for the field of rational numbers:

\begin{example} \label{EX:Q2}
For a finite subgroup $\G \subseteq \GL_n(\QQ)$ and a given prime
$\ell$, Proposition~\ref{P:p'} implies that the $\ell$-part
$|\G|_{\ell}$ of the order of $\G$ divides
$|\GL_n(\FF_{\!p})|_{\ell}$, where $p$ is any prime other than
$\ell$. Furthermore, if $\ell \neq 2$ then
$|\GL_n(\FF_{\!p})|_{\ell} = M(n)_\ell$ for infinitely many primes
$p$, by \eqref{E:p1}. Thus, we have shown (again) that if $\G$ is a
finite subgroup of $\GL_n(\QQ)$ then \text{$|\G|_{\ell}$ divides
$M(n)_{\ell}$ for all primes $\ell \neq 2$.} In order to extend this
to the prime $\ell = 2$, Minkowski uses additional facts about
quadratic forms. This will be explained below.
\end{example}

\subsection{The Schur bound} \label{SS:schurbound}

Fix an algebraic number field $K$. We will describe certain
constants $S(n,K)$, introduced by Schur in \cite{iS05}, for the
purpose of extending Theorem~\ref{T:schur} to general algebraic
number fields. Thus, $S(n,\QQ)$ will be identical to $M(n)$. Like
$M(n)$, the constant $S(n,K)$ will be defined as a product of
$\ell$-factors for all prime numbers $\ell$, and almost all
$\ell$-factors will be $1$. Throughout, we put
$$
\zeta_m = e^{2\pi i/m}\in \CC \ .
$$

For a given prime $\ell$, the chain $K \cap \QQ(\zeta_{\ell})
\subseteq \dots \subseteq K \cap \QQ(\zeta_{\ell^m}) \subseteq K
\cap \QQ(\zeta_{\ell^{m+1}}) \subseteq \dots$ of subfields of $K$
must stabilize, since $K$ is finite over $\QQ$. Thus we may define
\begin{equation} \label{E:ml}
m(K,\ell) = \min\{m \ge 1 \mid  K \cap \QQ(\zeta_{\ell^m}) = K \cap
\QQ(\zeta_{\ell^{m+1}}) = \dots \} \ .
\end{equation}
Now put
\begin{equation} \label{E:tl}
t(K,\ell) = [\QQ(\zeta_{\ell^{m(K,\ell)}}) : K \cap
\QQ(\zeta_{\ell^{m(K,\ell)}})]\ .
\end{equation}
and define
\begin{align}
S(n,K) &= 2^{n - \left\lfloor \frac{n}{t(K,2)}\right\rfloor}
\prod_{\ell} \ell^{m(K,\ell)\left\lfloor
\frac{n}{t(K,\ell)}\right\rfloor + \left\lfloor \frac{n}{\ell
t(K,\ell)}\right\rfloor + \left\lfloor \frac{n}{\ell^2
t(K,\ell)}\right\rfloor + \dots} \label{E:schurbound} \\
&= 2^{n - \left\lfloor \frac{n}{t(K,2)}\right\rfloor} \prod_{\ell}
\ell^{m(K,\ell)\left\lfloor \frac{n}{t(K,\ell)}\right\rfloor}
\left(\left\lfloor \tfrac{n}{t(K,\ell)}\right\rfloor
!\right)_{\ell}\ . \notag
\end{align}
Here, $\ell$ runs over all rational primes, including $2$, and the
second equality follows from equation \eqref{E:factorial}. Since
$t(K,\ell)[K:\QQ] \ge \ell -1$, only finitely many $\ell$ will
satisfy $t(K,\ell) \le n$ and so almost all $\ell$-factors are
trivial.

\begin{example} \label{EX:Q3}
Let $K = \QQ(\zeta_k)$ for some positive integer $k$. Since
$\QQ(\zeta_k) \cap \QQ(\zeta_t) = \QQ(\zeta_{(k,t)})$, we have
$m(K,\ell) = \max\{1, v_{\ell}(k)\}$. If $\ell$ does not divide $k$
then $t(K,\ell) = \ell-1$; otherwise $t(K,\ell)=1$. For $K = \QQ$ in
particular, we obtain $m(\QQ,\ell) = 1$ and $t(\QQ,\ell) = \ell - 1$
for all $\ell$. Thus, equation \eqref{E:schurbound} reduces to
\eqref{E:minkowskibound} and so $S(n,\QQ) = M(n)$.
\end{example}

In \cite{iS05}, Schur proved the following generalization of
Theorem~\ref{T:schur} using a larger dose of character theory than
what was needed in Section~\ref{S:proof}.

\begin{thm}[Schur 1905] \label{T:schur2}
Let $\G$ be a finite subgroup of $\GL_n(\CC)$ such that the traces
of all elements of $\G$ belong to some fixed algebraic number field
$K$. Then $|\G|$ divides $S(n,K)$.
\end{thm}

An alternative description of the constants $S(n,K)$ is as follows.
Let $\bmu_{\ell^\infty}$ denote the group of all $\ell$-power
complex roots of unity. Then $K \cap \QQ(\bmu_{\ell^\infty}) = K
\cap \QQ(\zeta_{\ell^{m(K,\ell)}})$.
\medskip

$\bullet$ If $\ell$ is odd then each $K \cap
\QQ(\zeta_{\ell^m})/\QQ$ is a subextension of
$\QQ(\zeta_{\ell^m})/\QQ$ which is cyclic with Galois group
isomorphic to $(\ZZ/\ell\ZZ)^* \cong \ZZ/\ell^{m-1}\ZZ \times
\ZZ/(\ell - 1)\ZZ$. Also, $K \cap \QQ(\zeta_{\ell})$ is the fixed
subfield of $K \cap \QQ(\zeta_{\ell^m})$ under the group
$\ZZ/\ell^{m-1}\ZZ$. Thus, $[K \cap \QQ(\zeta_{\ell^m}) : \QQ] = [K
\cap \QQ(\zeta_{\ell}) : \QQ][K \cap \QQ(\zeta_{\ell^m}) :
\QQ]_{\ell}$ and $[K \cap \QQ(\zeta_{\ell}) : \QQ]$ is a divisor of
$\ell-1$. Hence, for odd primes $\ell$,
\begin{align}
m(K,\ell) &= 1 + v_{\ell}([ K \cap \QQ(\bmu_{\ell^\infty}) : \QQ ])
\label{E:mtodd} \\
t(K,\ell) &= [ \QQ(\zeta_\ell) : K \cap \QQ(\zeta_\ell)] =
\tfrac{\ell - 1}{(\ell - 1,[ K \cap \QQ(\bmu_{\ell^\infty}) : \QQ
])}\notag \ .
\end{align}

$\bullet$ For the prime $\ell = 2$, the extension
$\QQ(\zeta_{2^m})/\QQ$ has Galois group $(\ZZ/2^m\ZZ)^* \cong
\ZZ/2^{m-2}\ZZ \times \ZZ/2\ZZ$  ($m \ge 2$). The factor $\ZZ/2\ZZ$
is generated by complex conjugation. When $m > 2$, the field
$\QQ(\zeta_{2^m})$ has exactly three subfields that are not
contained in $\QQ(\zeta_{2^{m-1}})$: besides $\QQ(\zeta_{2^m})$,
there are $\QQ(\zeta_{2^m} + \zeta_{2^m}^{-1})$ and $\QQ(\zeta_{2^m}
- \zeta_{2^m}^{-1})$. If $t(K,2) = 1$, which certainly holds when
$m(K,2) = 1$ or $m(K,2) = 2$, then $K \cap \QQ(\bmu_{2^\infty}) =
\QQ(\zeta_{2^{m(K,2)}})$, and so $[K \cap \QQ(\bmu_{2^\infty}) :
\QQ] = 2^{m(K,2)-1}$. If $t(K,2) \neq 1$ then $K \cap
\QQ(\bmu_{2^\infty})$ must be equal to either
$\QQ(\zeta_{2^{m(K,2)}} + \zeta_{2^{m(K,2)}}^{-1})$ or
$\QQ(\zeta_{2^{m(K,2)}} - \zeta_{2^{m(K,2)}}^{-1})$. Thus,
$t(K,2)=2$ and $[K \cap \QQ(\bmu_{2^\infty}) : \QQ] = 2^{m(K,2)-2}$.
In either case, the $2$-factor of $S(n,K)$ in \eqref{E:schurbound}
simplifies to
\begin{equation} \label{E:schurbound2}
S(n,K)_2 = [K \cap \QQ(\bmu_{2^\infty}) : \QQ]^{\left\lfloor
\frac{n}{t(K,2)}\right\rfloor} 2^n (n!)_2
\end{equation}

The following properties of $S(n,K)$ are easy to verify:
\begin{equation} \label{E:schurbound3}
S(m,K)S(n,K) \text{ divides } S(m+n,K)
\end{equation}
and
\begin{equation} \label{E:schurbound4}
\text{$S(n,K)$ divides $S(n,F)$ if $K \subseteq F$.}
\end{equation}

\subsection{Odd primes}
\label{SS:oddprimes}

The following proposition establishes Theorem~\ref{T:schur2} for the
$2'$-part of $|\G|$. The special case where $K = \QQ$ was done
earlier in Example~\ref{EX:Q2}.

\begin{prop} \label{P:oddell}
Let $\G$ be a finite subgroup of $\GL_n(\CC)$. Assume that the
traces of all elements of $\G$ belong to some algebraic number field
$K$. Then $|\G|_{\ell}$ divides $S(n,K)$ for all odd primes $\ell$.
\end{prop}

\begin{proof}
Replacing $\G$ by a conjugate in $\GL_n(\CC)$ if necessary, we can
make sure that $\G \subseteq \GL_n(F)$ for some algebraic number
field $F \supseteq K$. Indeed, any splitting field for $\G$ that is
finite over $K$ will serve this purpose; see \cite[Theorem
9.9]{mI76}. Let $\cO = \cO_{F}$ denote the ring of algebraic
integers of $F$ and consider any non-zero prime $\p$ of $\cO$ such
that $\G \subseteq \GL_n(\cO_{\p})$ and $\ell \notin \p$. Put $(p) =
\p \cap \ZZ$ and assume that $p$ is chosen as in Lemma~\ref{L:Q} and
also satisfies $p
> n$. Let $\rho \colon \G \to
\GL_n(\FF_{\!\N(\p)})$ denote the reduction homomorphism
\eqref{E:reduction} restricted to $\G$. Upon replacing $\G$ by a
Sylow $\ell$-subgroup, the map $\rho$ becomes injective, by
Lemma~\ref{L:p'}, and our goal now is to show that $|\G|$ divides
$S(n,K)_{\ell}$.

As in the first paragraph of the proof of Lemma~\ref{L:algebraic},
one sees that the traces of all elements of $\G$ actually belong to
the ring of algebraic integers $\cO_{K'}$ of the field $K' = K \cap
\QQ(\bmu_{\ell^\infty})$. Therefore, $\tr(\rho(g)) \in \FFq$ holds
for all $g \in \G$, where $q = \N(\p \cap \cO_{K'}) = p^{f}$.
Lemma~\ref{L:galois} now implies that $\rho(\G)$ is conjugate to a
subgroup of $\GL_n(\FFq)$ and Lemma~\ref{L:Q} further gives that
$$
|\G| \quad \text{divides} \quad |\GL_n(\FFq)|_{\ell} = \ell^{\left(
1 + v_{\ell}(f)\right)\left\lfloor \frac{n}{\tau} \right\rfloor}
\left(\left\lfloor \tfrac{n}{\tau} \right\rfloor !\right)_{\ell} \ ,
$$
where $\tau = \frac{\ell-1}{(\ell - 1,f)}$. Now, for odd $\ell$,
$$
S(n,K)_{\ell} = \ell^{m(K,\ell)\left\lfloor
\frac{n}{t(K,\ell)}\right\rfloor} \left(\left\lfloor
\tfrac{n}{t(K,\ell)}\right\rfloor !\right)_{\ell} \\
$$
with $m(K,\ell) = 1 + v_{\ell}([ K \cap \QQ(\bmu_{\ell^\infty}) :
\QQ ])$ and $t(K,\ell) = \frac{\ell - 1}{(\ell - 1,[ K \cap
\QQ(\bmu_{\ell^\infty}) : \QQ ])}$ by \eqref{E:mtodd}. Since the
residue class of $p$ generates $(\ZZ/\ell^{s}\ZZ)^*$ for all $s$,
$p$ remains prime in $\ZZ[\zeta_{\ell^{s}}]$; see the proof of
Lemma~\ref{L:Q} and \cite[Theorem 2 on p.~196]{kImR93}. In
particular, $p$ remains prime in $\cO_{K'}$, and so $f = f(\p \cap
\cO_{K'}/\QQ) = [ K \cap \QQ(\bmu_{\ell^\infty}) : \QQ ]$.
Therefore, $|\GL_n(\FFq)|_{\ell} = S(n,K)_{\ell}$ and the
proposition is proved.
\end{proof}

\subsection{Unitary, orthogonal and symplectic groups}
\label{SS:bilinear}

In this section, we review some standard facts about hermitian and
skew-hermitian forms and certain classical groups that are
associated with them. Throughout, $\k$ will denote a field and $a
\mapsto a^\theta$ will be an automorphism of $\k$ satisfying
$\theta^2 = \Id$. We assume for simplicity that $\ch \k \neq 2$.

\subsubsection{Sesquilinear forms} \label{SSS:1.5}

Let $V$ denote an $n$-dimensional vector space over $\k$. A
bi-additive map $\beta \colon V \times V \to \k$ is called
\emph{sesquilinear} (with respect to $\theta$) if
$$
\beta(av,bw) = ab^{\theta}\beta(v,w)
$$
holds for all $v,w\in V$ and $a,b \in \k$. When $\theta$ is the
identity, sesquilinear forms are ordinary bilinear forms. A
sesquilinear form $\beta$ is called \emph{non-singular} if $\beta$
satisfies the following equivalent conditions: (i) $\beta(v,V) =
\{0\}$ for $v \in V$ implies $v=0$; (ii) $\beta(V,v) = \{0\}$ for $v
\in V$ implies $v=0$; (iii) for any basis $\{v_1,\dots,v_n\}$ of
$V$, the matrix $\left(\beta(v_i,v_j)\right)_{n \times n}$ has
non-zero determinant; see \cite[Proposition XIII.7.2]{sL02}. If
$\beta$ is any sesquilinear form on $V$ and $g \in \GL(V)$ then,
defining $\beta^g(v,v') := \beta(g(v),g(v'))$ for $v,v' \in V$, one
again obtains a sesquilinear form $\beta^g$ on $V$ with respect to
$\theta$; it is called \emph{equivalent} to $\beta$.

Sesquilinear forms $\beta$ satisfying $\beta(w,v) =
\beta(v,w)^\theta$ (resp. $\beta(w,v) = -\beta(v,w)^\theta$) for all
$v,w\in V$ are called \emph{hermitian} (resp.
\emph{skew-hermitian}). The stabilizer in $\GL(V)$ of a non-singular
hermitian or skew-hermitian form $\beta$ is called the group of
\emph{isometries} of $(V,\beta)$ and is denoted by $\Iso(V,\beta)$;
so
$$
\Iso(V,\beta) = \{ g \in \GL(V) \mid \beta(g(v),g(v')) = \beta(v,v')
\text{ for all $v,v' \in V$}\} \ .
$$
Let $\beta$ be non-singular skew-hermitian. If $\beta(v,v) \neq 0$
for some $v \in V$ then $\beta' = \beta(v,v)\beta$ is a non-singular
hermitian form on $V$ with $\Iso(V,\beta') = \Iso(V,\beta)$. On the
other hand, if $\beta(v,v) = 0$ for all $v \in V$ then it is easy to
see that $\theta = \Id$ and so $\beta$ is an alternating bilinear
form. Therefore, when studying isometry groups of non-singular
hermitian or skew-hermitian forms $\beta$ on $V$, it suffices to
consider the following cases:
\begin{description}
\item[unitary case] $\beta$ is hermitian with respect to $\theta
\neq \Id$;
\item[orthogonal case] $\beta$ is symmetric bilinear ($\theta =
\Id$);
\item[symplectic case] $\beta$ is alternating bilinear ($\theta =
\Id$).
\end{description}

\subsubsection{Twisting modules} \label{SSS:twisting}

Now assume that $V$ is a finitely generated (left) $\k[\G]$-module,
where $\G$ is a finite group. We let $V^\theta = \{ v^\theta \mid v
\in V\}$ denote a copy of $V$  with operations
$$
v^\theta + w^\theta = (v+w)^\theta \,,\ \ (av)^\theta = a^\theta
v^\theta \quad \text{and} \quad g v^\theta = (gv)^\theta
$$
for $v,w \in V$, $a \in \k$ and $g \in \G$. Then $V^\theta$ becomes
a $\k[\G]$-module and
\begin{equation} \label{E:twistedtrace}
\tr_{V^\theta/\k}(g) = \left( \tr_{V/\k}(g) \right)^\theta
\end{equation}
holds for all $g \in \G$. Furthermore, there is an isomorphism of
$\k[\G]$-modules
\begin{equation} \label{E:twistediso}
\left(V \otimes_{\k} V^{\theta} \right)^* \cong \left\{
\text{sesquilinear forms $V \times V \to \k$ with respect to
$\theta$} \right\} \ .
\end{equation}
The isomorphism sends a linear for $\varphi \colon V \otimes_{\k}
V^{\theta} \to \k$ to the form $\widetilde{\varphi} \colon V \times
V \to \k$ given by $\widetilde{\varphi}(v,w) = \varphi(v \otimes
w^\theta)$. The group $\Sy_2 = \langle \tau \rangle$ acts on the
space of sesquilinear forms $\beta \colon V \times V \to \k$ with
respect to $\theta$ by
$$
(\tau \beta)(v,w) = \beta(v,w)^\theta
$$
for $v,w \in V$. This action commutes with the action of $\G$. Note
however that the action is only $\k$-semilinear: $\tau(a\beta) =
a^{\theta}\tau\beta$. Clearly, $\beta$ is hermitian (resp.
skew-hermitian) if and only if $\tau\beta = \beta$ (resp. $\tau\beta
= -\beta$).

\begin{lem} \label{L:twisting}
Let $\sigma \colon \G \to \GL(V)$ be an irreducible representation
of the finite group $\G$. If $V^* \cong V^\theta$ as
$\k[\G]$-modules then $\sigma(\G) \subseteq \Iso(V,\beta)$ for some
non-singular form $\beta$ on $V$ that is hermitian or skew-hermitian
with respect to $\theta$.
\end{lem}

\begin{proof}
Since $V^* \cong V^\theta$, we have $V^* \otimes_{\k} V \cong \left(V
\otimes_{\k} V^{\theta} \right)^*$ and so
$$
\End_{\k}(V) \cong \left\{ \text{sesquilinear forms $V \times V \to
\k$ with respect to $\theta$} \right\}
$$
as  $\k[\G]$-modules, by \eqref{E:twistediso}. The identity $\Id_V
\in \End_{\k}(V)$ therefore corresponds to a non-zero $\G$-invariant
sesquilinear form $\beta$. Write $\beta = \beta_+ + \beta_-$ with
$\beta_\pm = \frac{1}{2}(1 \pm \tau)(\beta)$, where $\Sy_2 = \langle
\tau \rangle$ as above. Then $\tau\beta_\pm = \pm\beta$; so
$\beta_+$ is hermitian and $\beta_-$ is skew-hermitian with respect
to $\theta$, and at least one of them is non-zero. Moreover, both
$\beta_\pm$ are $\G$-invariant, since the actions of $\tau$ and $\G$
commute. Finally, any non-zero $\G$-invariant hermitian or
skew-hermitian form on $V$ is non-singular, because its radical is a
proper $\k[\G]$-submodule of $V$, and hence it must be zero because
$V$ is assumed simple.
\end{proof}

\subsubsection{Isometry groups over finite fields} \label{SSS:isometry}

We will now concentrate on the case of a finite field $\k = \FFq$ of
order $q = p^f$ for some odd prime $p$. Let $\beta$ be a
non-singular hermitian or skew-hermitian form on $V \cong \FFq^n$.
Since we are only interested in the group of isometries
$\Iso(V,\beta)$, we may assume that $\beta$ is unitary, orthogonal
or symplectic. The orders of these groups are classical; see
Dieudonn{\'e} \cite{jDi71} or Artin \cite[Section III.6]{eA88}, for
example. The original sources are Minkowski's dissertation
\cite{hM85} and Dickson \cite{leD01}.

\begin{description}
\item[unitary case] Since $\theta$ has order $2$ in this case, $f$ must
be even. Moreover, $\beta$ is unique up to equivalence, and so
$\Iso(V,\beta)$ is determined up to conjugation. The order of
$\Iso(V,\beta)$ is
\begin{equation} \label{E:Uorder}
|\Iso(V,\beta)| = p^{fn(n-1)/4} \prod_{i=1}^n (p^{fi/2} - (-1)^i) \
.
\end{equation}
\item[symplectic case] Again, $\beta$ is unique up to
equivalence. The dimension $n$ must be even. One has
\begin{equation} \label{E:SPorder}
|\Iso(V,\beta)| = q^{n^2/4} \prod_{i=1}^{n/2} (q^{2i} - 1)\ .
\end{equation}
\item[orthogonal case] Here, the order of $\Iso(V,\beta)$ is given by
\begin{equation} \label{E:Oorder}
|\Iso(V,\beta)| = \begin{cases}
\ 2q^{ (n-1)^2/4} \displaystyle\prod_{i=1}^{(n-1)/2} (q^{2i} - 1) & \text{if $n$ is odd,} \\
\ 2q^{n(n-2)/4}(q^{n/2} - \varepsilon)
\displaystyle\prod_{i=1}^{(n-2)/2} (q^{2i} - 1) & \text{if $n$ is
even,}
\end{cases}
\end{equation}
where $\varepsilon = \pm 1$ depends on the form $\beta$. The
detailed description of $\varepsilon$ will not matter for us.
\end{description}

\begin{lem} \label{L:isoestimate}
Let $K$ be an algebraic number field contained in
$\QQ(\bmu_{2^\infty})$ (so $K$ is Galois over $\QQ$ and in
particular stable under complex conjugation). If $K \nsubseteq \RR$
then assume that $t(K,2)=1$. There are infinitely many odd primes
$\p$ of the ring of algebraic integers $\cO_{K}$ satisfying the
following two conditions:
\begin{enumerate}
\item[(i)] $\p$ is stable under complex conjugation, and
\item[(ii)] If $\beta$ is any non-singular
hermitian or skew-hermitian form on $V = \FFq^n$ with respect to the
automorphism $\theta$ of $\cO_{K}/\p = \FF_q$ that is afforded by
complex conjugation then $|\Iso(V,\beta)|_2$ divides $S(n,K)_2$.
\end{enumerate}
\end{lem}

\begin{proof}
We will need the following elementary observation. If $p$ is a prime
satisfying $p \equiv -1 + 2^k \mod 2^{k+1}$ for some $k \ge 2$ then,
for all positive integers $i$,
\begin{equation} \label{E:p2i-1}
(p^{i} - (-1)^i)_2 = 2^{k}i_2\ .
\end{equation}
To see this, we remark first that $(p^i-1)_2 = 2$ holds for odd$i$,
because the residue class of $p$ modulo $4$ is the nonidentity
element of $(\ZZ/4\ZZ)^*$, and hence the same holds for all odd
powers of $p$. Moreover, since $p^{2} \equiv 1 \mod 2^{k+1}$, we
have $p^{i} \equiv 1 \mod 2^{k+1}$ for all even $i$, and hence
$(p^i+1)_2 = 2$. Now, to prove \eqref{E:p2i-1}, assume first that
$i$ is odd, say $i = 2j+1$. Then the foregoing implies that $p^{i} -
(-1)^i = p^{2j}p + 1 \equiv p + 1 \equiv  2^k \mod 2^{k+1}$, and so
$(p^{i} - (-1)^i)_2 = 2^k$, proving \eqref{E:p2i-1} for odd values
of $i$. Finally, assume that $i=2j$. Then $p^{i} - (-1)^i =
(p^j-1)(p^j+1)$. If $j$ is odd then we know that $(p^j+1)_2 = 2^k$
and $(p^j-1)_2 = 2$, and hence $(p^{i} - (-1)^i)_2 = 2^{k+1}$, as
desired. When $j$ is even then $(p^j-1)_2 = 2^k j_2$, by induction,
and $(p^j+1)_2 = 2$, as we remarked earlier. Thus, \eqref{E:p2i-1}
is proved in all cases.

Turning to the proof of the lemma, note that $K/\QQ$ is Galois,
being a subextension of the abelian extension
$\QQ(\bmu_{2^\infty})/\QQ$. Put $m = m(K,2)$, $t = t(K,2)$ and
$\zeta = \zeta_{2^m}$. Then \eqref{E:schurbound2} becomes
$$
S(n,K)_2 = [K : \QQ]^{\left\lfloor n/t \right\rfloor} 2^n (n!)_2
$$
and $K$ is one of the fields $\QQ(\zeta)$ or $\QQ(\zeta +
\zeta^{-1})$; see \S\ref{SS:schurbound}. We will deal with each of
these cases separately. Throughout, $\p$ will denote a prime ideal
of $\cO_{K}$ and we put $q = \N(\p)$ and $(p) = \p \cap \ZZ$.

First consider the case where $K$ is real. Then property (i) is
automatic and $\Iso(V,\beta)$ is symplectic or orthogonal. Replacing
the factor $(q^{n/2} - \varepsilon)$ in formula \eqref{E:Oorder} for
even $n$ by its multiple $(q^{n/2} - \varepsilon)(q^{n/2} +
\varepsilon)/2  = (q^{n} - 1)/2$ and deleting $q$-factors (which are
odd) we obtain the expression $\prod_{i=1}^{n/2} (q^{2i} - 1)$ that
only depends on $n$ and $q$ and is identical to \eqref{E:SPorder}
stripped of its $q$-factors. Put
\begin{equation*} \label{E:on3}
\mathbf{o}(n,q) = \begin{cases}
\quad 2\displaystyle\prod_{i=1}^{(n-1)/2} (q^{2i} - 1) & \text{if $n$ is odd,} \\
\quad \displaystyle\prod_{i=1}^{n/2} (q^{2i} - 1) & \text{if $n$ is
even.}
\end{cases}
\end{equation*}
Now $q = p^f$, where $f = [\cO_{K}/\p :\FF_{\!p}]$ is a divisor of
$[K:\QQ]$; so $f$ is a power of $2$. Choose $\p$ to lie over any
rational prime $p$ with $p \equiv 3 \mod 8$. Then \eqref{E:p2i-1}
with $k=2$ implies that the $2$-part of $q^{2i} - 1$ for $i \ge 1$
is given by $\left(q^{2i} - 1\right)_2 = 8f i_2$. It follows that
the $2$-part of $\mathbf{o}(n,q)$ can be written as
$\mathbf{o}(n,q)_2 = f^{\left\lfloor n/2 \right\rfloor} 2^n (n!)_2$
in both cases. Since $f$ is a divisor of $[K:\QQ]$ and $t$ equals
$1$ or $2$, we see that $\mathbf{o}(n,q)_2$ divides $S(n,K)_2$ which
settles the symplectic and orthogonal cases.

Next, let $K = \QQ(\zeta)$ with $m \ge 2$. Choose $\p$ to lie over
any rational prime $p$ satisfying $p \equiv -1 + 2^m \mod 2^{m+1}$.
The $\p$ is stable under complex conjugation. Indeed, the
decomposition group of $\p$ is generated by the automorphism of $K$
sending $\zeta$ to $\zeta^p$ (cf., e.g.,  \cite[Corollary on
p.~197]{kImR93}), and our choice of $p$ implies that $\zeta^p =
\zeta^{-1} = \overline{\zeta}$. Thus, complex conjugation
$\overline{\phantom{I}}$ belongs to the decomposition group of $\p$,
and it must in fact generated the decomposition group, because
$\overline{\phantom{I}}$ is not a square in $\Gal(K/\QQ)$. Since
$\p$ is unramified over $\QQ$, its relative degree over $\QQ$ equals
$f=2$; so $q = p^2$. Therefore, \eqref{E:Uorder} and \eqref{E:p2i-1}
give
$$
|\Iso(V,\beta)|_2 =  \prod_{i=1}^n (p^{i} - (-1)^i)_2  = 2^{mn}
(n!)_2 = 2^{(m-1)n} 2^n (n!)_2\ .
$$
Since $[K:\QQ] = 2^{m-1}$ and $t=1$, the last expression is equal to
$S(n,K)_2$, thereby completing the proof of the lemma.
\end{proof}

The lemma fails in the excluded case $K \nsubseteq \RR$, $t(K,2)=2$.
For example, let $K = \QQ(\sqrt{-2})$. Then $m(K,2) = 3$ and $t(K,2)
= 2$ and so $S(n,K)_2 = 2^{\left\lfloor \frac{n}{2}\right\rfloor}
2^n (n!)_2$. On the other hand, if $\p$ is an odd prime of $\cO_K$
that is stable under complex conjugation, then $f(\p/\QQ) = 2$ and
$p \equiv -1 \mod 8$. It follows that $|\Iso(V,\beta)|_2 =
\prod_{i=1}^n (p^{i} - (-1)^i)_2$ is divisible by $2^{3n}$ which is
too big.

\subsection{The prime $\ell = 2$}
\label{SS:congruence3}

The following proposition complements Proposition~\ref{P:oddell}. It
would be nice to remove the restrictions $K' = K \cap
\QQ(\bmu_{2^\infty}) \subseteq \RR$ or $t(K,2) = 1$ on $K$. This
would require replacing the isometry groups $\Iso(V,\beta)$ by
suitable subgroups.

\begin{prop} \label{P:schur2}
Let $\G$ be a finite subgroup of $\GL_n(\CC)$ such that the traces
of all elements of $\G$ belong to some fixed algebraic number field
$K$. Assume that $K' = K \cap \QQ(\bmu_{2^\infty}) \subseteq \RR$ or
$t(K,2) = 1$. Then $|\G|_2$ divides $S(n,K)$.
\end{prop}

\begin{proof}
We may assume that $\G$ is a $2$-group. Therefore, $\tr(g) \in
\cO_{K'}$ for all $g \in \G$. Replacing $K$ by $K'$, we may assume
that $K = K' \subseteq \QQ(\bmu_{2^\infty})$; see
\eqref{E:schurbound4}. Choose a prime $\p$ of $\cO_K$ as in
Lemma~\ref{L:isoestimate} and put $q = \N(\p)$.  As in the proof of
Proposition~\ref{P:oddell}, we can arrange that $\G \subseteq
\GL_n(F)$ for some algebraic number field $F$ containing $K$. Choose
a prime $\bP$ of $\cO = \cO_F$ lying over $\p$ and put $\k =
\cO/\bP$; so $\FF_q = \cO_K/\p \subseteq \k$. We may assume $\G
\subseteq \GL_n(\cO_{\bP})$ and that $(p) = \bP \cap \ZZ$ satisfies
$p > n$. By Lemma~\ref{L:p'}, the reduction homomorphism
$\GL_n(\cO_{\bP}) \to \GL_n(\k)$ is injective on $\G$. We will write
the restriction of this map to $\G$ as
$$
\rho \colon \G \hookrightarrow \GL_n(\k) \ .
$$
Then $\tr \rho(g) = \tr g \mod \p \in \FF_q \subseteq \k$ for $g \in
\G$, and $ \tr \rho(g^{-1}) = \left(\tr \rho(g) \right)^{\theta}$,
where $\theta$ denotes the automorphism of $\FF_q$ that is afforded
by complex conjugation, as in Lemma~\ref{L:isoestimate}. Now
Lemma~\ref{L:galois} implies that $\rho(\G)^v = v^{-1}\rho(\G)v
\subseteq \GL_n(\FFq)$ for some $v \in \GL_n(\k)$; so we may
consider the representation
$$
\sigma = (\,.\,)^v \circ \rho  \colon \G \hookrightarrow \GL(V) \ ,
$$
where $V = \FFq^n$. Note that $\tr \sigma(g) = \tr \rho(g)$ for all
$g \in \G$. We will write $V$ as a direct sum of
$\FFq[\G]$-submodules $U_i$ on which $\G$ acts as a subgroup of
$\Iso(U_i,\beta_i)$ for some non-singular hermitian or
skew-hermitian form $\beta_i$ with respect to $\theta$ on $U_i$.
This will imply that $|\G|$ divides $\prod_i |\Iso(U_i,\beta_i)|_2$,
and hence $|\G|$ divides $\prod_i S(\dim U_i,K)$ by
Lemma~\ref{L:isoestimate}. Since $\prod_i S(\dim U_i,K)$ is a
divisor of $S(\sum_i \dim U_i,K) = S(n,K)$, by
\eqref{E:schurbound3}, the theorem will follow.

To achieve the decomposition of $V$, recall that $\tr \sigma(g^{-1})
= \left(\tr \sigma(g) \right)^{\theta}$ for all $g \in \G$. By
\eqref{E:twistedtrace}, this says that the $\FFq[\G]$-modules $V^*$
and $V^\theta$ have the same character, and hence they are
isomorphic; see the proof of Lemma~\ref{L:galois}. Write $V \cong
\bigoplus_i V_i^{(n_i)}$ with non-isomorphic irreducible
$\FFq[\G]$-modules $V_i$. Then $V^* \cong \bigoplus_i
\left(V_i^*\right)^{(n_i)}$ and $V^\theta \cong \bigoplus_i
\left(V_i^{\theta}\right)^{(n_i)}$. For each $i$, there is an $i'$
so that $V_i^* \cong V_{i'}^\theta$. If $i = i'$ then
Lemma~\ref{L:twisting} says that $\G$ acts on $V_i$ as a subgroup of
$\Iso(V_i,\beta_i)$ for some non-singular hermitian or
skew-hermitian form $\beta_i$ on $V_i$. Now assume that $i \neq i'$.
Then $V_i^* \oplus V_i^\theta$ is a direct summand of $V^\theta$,
and hence $\widetilde{V_i} = \left(V_i^*\right)^{\theta} \oplus V_i$
is a direct summand of $V$. Defining
$$
\beta_i(f^{\theta} + v, f'^{\theta} + v') = f(v')^{\theta} + f'(v)
$$
for $f,f' \in V_i^*$ and $v,v' \in V_i$ we obtain a non-singular
hermitian form on $\widetilde{V_i}$ that is preserved by the action
of $\G$. This yields the desired decomposition of $V$ and completes
the proof of the theorem.
\end{proof}

\newpage

\section{Outlook} \label{S:outlook}

We conclude by surveying, without proofs, a number of topics that
are related to the foregoing.

\subsection{The largest groups and recent work on the Jordan bound}
\label{SS:largest}

\subsubsection{}

The group $\G$ constructed in Proposition~\ref{P:minkowski} is
isomorphic to the so-called wreath product
$$
\Sy_{m+1} \wr \Sy_a \ .
$$
By definition, $\Sy_{m+1} \wr \Sy_a$ is the semidirect product of
$\Sy_{m+1}^a \rtimes \Sy_{a}$\,, where $\Sy_a$ acts on $\Sy_{m+1}^a
= \Sy_{m+1} \times \dots \times \Sy_{m+1}$ by permuting the $a$
factors $\Sy_{m+1}$. The special case $m=1$ yields the group $\{\pm
1\} \wr \Sy_n$, a subgroup of $\GL_n(\ZZ)$ order $2^n n!$ which is
also known as the automorphism group $\Aut(B_n)$ of the root system
of type $B_n$; see \cite{nB68}. For almost all values of $n$, these
particular groups turn out to be the largest finite groups that can
be found inside $\GL_n(\ZZ)$, and even inside $\GL_n(\QQ)$ (see
\S\ref{SS:congruence}). Indeed, Feit \cite{wFxx} has shown that, for
all $n > 10$ and for $n = 1,3,5$, the finite subgroups of
$\GL_n(\QQ)$ of largest order are precisely the conjugates of
$\Aut(B_n)$. For the remaining values of $n$, Feit also
characterizes the largest finite subgroups of $\GL_n(\QQ)$ and shows
that they are unique up to conjugacy. Feit's proof depends in an
essential way on an unfinished manuscript of Weisfeiler \cite{bWxx}
which establishes an estimate for the so-called Jordan bound; see
\S~\ref{SSS:jordan} below. An alternative proof of Feit's theorem
for sufficiently large values of $n$ has been given by Friedland
\cite{sF97} who relies on another (published)  article of
Weisfeiler's, \cite{bW84}. Both \cite{bWxx} and \cite{bW84} depend
crucially on the classification of finite simple groups.

Sadly, the two protagonists of the developments sketched above are
no longer with us: Walter Feit passed away on July 29, 2004 while
Boris Weisfeiler disappeared in January 1985 during a hiking trip in
the Chilean Andes. The present status of the investigation into
Weisfeiler's disappearance is documented on the web site
\textsl{http://www.weisfeiler.com/boris/}. For further information
on the subject of finite subgroups of $\GL_n(\ZZ)$ and of
$\GL_n(\QQ)$, especially maximal ones, see, e.g., Nebe and Plesken
\cite{gNwP95}, Plesken \cite{wP91}, the first chapter of \cite{mL05}
and, at a more elementary level, the article \cite{jKaP02} by
Kuzmanovich and Pavlichenkov.

\subsubsection{} \label{SSS:jordan}

The Jordan bound comes from the following classical result
\cite{cJ78}.

\begin{thm}[Jordan 1878] \label{T:jordan}
There exists a function $j \colon \NN \to \NN$  such that every
finite subgroup of $\GL_n(\CC)$ contains an abelian normal subgroup
of index at most $j(n)$.
\end{thm}

Early estimates for the optimal function $j(n)$ were quite
astronomical. Until fairly recently, the best known result was due
to Blichfeldt: $j(n) \le n!\, 6^{(n-1)(\pi(n+1)+1)}$, where
$\pi(n+1)$ denotes the number of primes $\le n+1$; see \cite[Theorem
30.4]{lD71}. Since $\Sy_{n+1} \into \GL_n(\CC)$, as explained in the
proof of Proposition~\ref{P:minkowski}, one must certainly have
$j(n) \ge (n+1)!$ for $n \ge 4$\,. In his near-complete manuscript
\cite{bWxx}, Weisfeiler comes close to proving that equality holds
for large enough $n$: he shows that if $n > 63$ then $j(n) \le
(n+2)!$. In \cite{bW84}, Weisfeiler announces the weaker upper bound
$j(n) \le n^{a \log n + b} n!$. Quite recently, Michael Collins
\cite{mC05b} was able to settle the problem by showing that for $n
\ge 71$ we do indeed have $j(n) = (n+1)!$ and, if this bound is
achieved by $\G$, then $\G$ modulo its center is isomorphic to
$\Sy_{n+1}$.

\subsubsection{} \label{SSS:jordanp}

Analogs of Jordan's Theorem for linear groups in characteristics $p
> 0$ were established by Weisfeiler \cite{bWxx}, \cite{bW84}, Larsen
and Pink \cite{mLrP98}, and Collins \cite{mC05c}. While both
Weisfeiler and Collins rely on the classification theorem, Larsen
and Pink prove a noneffective version of Jordan's theorem, without
explicit index and degree bounds, by using methods from algebraic
geometry and the theory of linear algebraic groups instead. We will
explain Collins' modular version of Jordan's Theorem. As usual,
$O_p(\G)$ denotes the maximal normal $p$-subgroup of the finite
group $\G$. Furthermore, a group is called \emph{quasisimple} if it
is perfect and simple modulo its center. Collins' result then reads
as follows.

\begin{thm}[Collins 2005] \label{T:jordanp}
Let $F$ be a field of positive characteristic $p$ and let $\G$ be a
finite subgroup of $\GL_n(F)$, where $n \ge 71$. Put $G =
\G/O_p(\G)$. Then $G$ has a normal subgroup $N$ such that
\begin{enumerate}
\item
$N = A\,Q_1 \dots Q_m$, a central product with $A$ abelian and the
$Q_i$ (quasi)simple Chevalley groups in characteristic $p$.
\item
$ [G : N] \le
\begin{cases}
\quad (n+2)! & \text{if $p$ divides $n+2$,} \\
\quad (n+1)! & \text{otherwise.}
\end{cases}
$
\end{enumerate}
\end{thm}

\subsection{The Minkowski sequence $M(n)$} \label{SS:sequence}

A search of Sloane's \emph{On-Line Encyclopedia of Integer
Sequences} \cite{OEIS}, by entering the first six terms
$2,24,48,5760,11520,2903040$ of $M(n)$, turns up a sequence labeled
A053657. This sequence has two additional descriptions besides
Minkowski's description of $M(n)$ as the least common multiple of
the orders of all finite subgroups of $\GL_n(\QQ)$; the other two
will be given below. We know of no direct argument explaining the
(proven) equivalence of $M(n)$ to the first sequence below. The
equivalence of the second sequence to $M(n)$ is currently supported
only by empirical evidence.

\begin{itemize}
\item
By Chabert et.~al.~\cite{CjlEtAl97}, the collection of all leading
coefficients of polynomials $f(x) \in \QQ[x]$ of degree at most $n$
such that $f(p) \in \ZZ$ holds for all primes $p$ is a fractional
ideal of the form $\frac{1}{a(n)}\ZZ$ for suitable positive integers
$a(n)$. It turns out that formula \eqref{E:minkowskibound} is
identical with the formula given in \cite[Proposition
4.1]{CjlEtAl97} for $a(n+1)$. Thus $M(n-1) = a(n)$.

\item
Following Paul Hanna \cite[A075264]{OEIS}, we let $P(n,z)$ denote
the coefficient of $x^n$ in the Taylor series for
$(\frac{-\ln(1-x)}{x})^z$ at $x=0$. Thus, $\sum_{m=1}^\infty
\binom{z}{m} \xi^m = \sum_{n=1}^\infty P(n,z)x^n$ with $\xi =
\frac{-\ln(1-x)}{x} - 1 = \sum_{k=1}^\infty \frac{x^k}{k+1}$ and
$\binom{z}{m} = \frac{z(z-1)\dots(z-m+1)}{m!}$\,. For example,
$P(1,z) = \frac{z}{2}$, $P(2,z) = \frac{5z + 3z^2}{24}$, $P(3,z) =
\frac{6z + 5z^2 + z^3}{48}$. In general, $P(n,z) \in z\QQ[z]$; the
polynomials $P(n,z)$ for $n \le 8$ are listed in sequence A075264 of
OEIS \cite{OEIS}. Paul Hanna has noted that the denominator of
$P(n,z)$, that is, the positive generator of the ideal $\{ q \in\ZZ
\mid q P(n,z) \in \ZZ[z]\}$, appears to coincide with $M(n)$.
\end{itemize}

In \cite{hM87}, Minkowski states the following recursion for the
sequence $M(n)$; the recursion is easy to check from
\eqref{E:minkowskibound}:
\begin{equation} \label{E:mn2}
M(2n+1) = 2\,M(2n) \qquad \text{and}\qquad M(2n) =
2\,M(2n-1)\,\prod_{p \colon p-1\mid 2n} p n_p \ .
\end{equation}
The product in \eqref{E:mn2} ranges over all primes $p$ such that
$p-1$ divides $2n$, and $n_p$ denotes the $p$-part of $n$, as usual.
This product has an interpretation in terms of the familiar
\emph{Bernoulli numbers} $B_n$ which are defined by $\frac{x}{e^x-1}
= \sum_{n=0}^\infty B_n\frac{x^n}{n!}$. In fact, $B_n = 0$ for odd
$n>1$ while $B_{2n}$ is a rational number whose denominator, when
written in lowest terms, is given by the von Staudt-Clausen theorem:
it is equal to $\prod_{p \colon p-1\mid 2n} p$\,; cf.~\cite[Theorem
1]{lC68}. Moreover, for each prime $p$ such that $p-1$ does not
divide $2n$, the numerator of $B_{2n}$ is divisible by the $p$-part
$n_p$\,; see~\cite[Theorem 5]{lC68}. Consequently, the product
$\prod_{p \colon p-1\mid 2n} p n_p$ in \eqref{E:mn2} is equal to the
denominator of $\frac{B_{2n}}{n}$\,. This was already pointed out by
Minkowski in \cite{hM87}. Finally, the asymptotic order of $M_n$ has
been determined by Katznelson \cite{yK94}: $\lim_{n \to \infty}
\left(M(n)/n!\right)^{1/n} = \prod_p p^{1/(p-1)^2} \approx 3.4109$.


\begin{ack}
The authors wish to thank Boris Datskovsky and Ed Letzter for their
comments on a preliminary version of this article, and Michael
Collins for making drafts of \cite{mC05a},\cite{mC05c},\cite{mC05b}
available to them.
\end{ack}


\bibliographystyle{amsplain}
\bibliography{bibliography}

\providecommand{\bysame}{\leavevmode\hbox to3em{\hrulefill}\thinspace}
\providecommand{\MR}{\relax\ifhmode\unskip\space\fi MR }
\providecommand{\MRhref}[2]{%
  \href{http://www.ams.org/mathscinet-getitem?mr=#1}{#2}
}
\providecommand{\href}[2]{#2}
\begin{thebibliography}{10}

\bibitem{eA88}
Emil Artin, \emph{Geometric algebra}, Wiley Classics Library, John Wiley \&
  Sons Inc., New York, 1988, Reprint of the 1957 original, A Wiley-Interscience
  Publication. \MR{90h:51003}

\bibitem{aB91}
Armand Borel, \emph{Linear algebraic groups}, second ed., Graduate Texts in
  Mathematics, vol. 126, Springer-Verlag, New York, 1991. \MR{92d:20001}

\bibitem{nB68}
Nicolas Bourbaki, \emph{{G}roupes et alg\`ebres de {L}ie. {C}hapitre {IV}:
  {G}roupes de {C}oxeter et syst\`emes de {T}its. {C}hapitre {V}: {G}roupes
  engendr\'es par des r\'eflexions. {C}hapitre {VI}: syst\`emes de racines},
  Actualit\'es Scientifiques et Industrielles, No. 1337, Hermann, Paris, 1968.
  \MR{39 \#1590}

\bibitem{nB73}
\bysame, \emph{\'{E}l\'ements de math\'ematique, {F}asc. {XXIII}}, Hermann,
  Paris, 1973, Livre II: Alg\`ebre. Chapitre 8: Modules et anneaux
  semi-simples, Nouveau tirage de l'\'edition de 1958, Actualit\'es
  Scientifiques et Industrielles, No. 1261. \MR{54 \#5282}

\bibitem{nB81}
\bysame, \emph{Alg\`ebre, chapitres 4 \`a 7}, Masson, Paris, 1981.
  \MR{84d:00002}

\bibitem{lC68}
Leonard Carlitz, \emph{Bernoulli numbers}, Fibonacci Quart \textbf{6} (1968),
  no.~3, 71--85. \MR{38 \#1071}

\bibitem{CjlEtAl97}
Jean-Luc Chabert, Scott~T. Chapman, and William~W. Smith, \emph{A basis for the
  ring of polynomials integer-valued on prime numbers}, Factorization in
  integral domains (Iowa City, IA, 1996), Lecture Notes in Pure and Appl.
  Math., vol. 189, Dekker, New York, 1997, pp.~271--284. \MR{99d:13008}

\bibitem{mC05a}
Michael Collins, \emph{Bounds for finite primitive complex linear groups},
  preprint, University of Oxford, 30 pages, 2005.

\bibitem{mC05c}
\bysame, \emph{Modular analogues of {J}ordan's theorem for finite linear
  groups}, preprint, University of Oxford, 43 pages, 2005.

\bibitem{mC05b}
\bysame, \emph{On {J}ordan's theorem for complex linear groups}, preprint,
  University of Oxford, 18 pages, 2005.

\bibitem{cCiR81}
Charles~W. Curtis and Irving Reiner, \emph{Methods of representation theory.
  {V}ol. {I}}, John Wiley \& Sons Inc., New York, 1981, With applications to
  finite groups and orders, Pure and Applied Mathematics, A Wiley-Interscience
  Publication. \MR{82i:20001}

\bibitem{leD01}
Leonard~Eugene Dickson, \emph{Linear groups: {W}ith an exposition of the
  {G}alois field theory}, Dover Publications Inc., New York, 1958, unaltered
  republication of the first edition [Teubner, Leipzig, 1901] with a new
  introduction by Wilhelm Magnus. \MR{21 \#3488}

\bibitem{jDi71}
Jean~A. Dieudonn{\'e}, \emph{La g\'eom\'etrie des groupes classiques},
  Springer-Verlag, Berlin, 1971, Troisi\`eme \'edition, Ergebnisse der
  Mathematik und ihrer Grenzgebiete, Band 5. \MR{46 \#9186}

\bibitem{lD71}
Larry Dornhoff, \emph{Group representation theory. {P}art {A}: {O}rdinary
  representation theory}, Marcel Dekker Inc., New York, 1971, Pure and Applied
  Mathematics, 7. \MR{50 \#458a}

\bibitem{wFxx}
Walter Feit, \emph{Orders of finite linear groups}, unpublished preprint,
  approx. 1998.

\bibitem{wF97}
Walter Feit, \emph{Finite linear groups and theorems of {M}inkowski and
  {S}chur}, Proc. Amer. Math. Soc. \textbf{125} (1997), no.~5, 1259--1262.
  \MR{97g:20007}

\bibitem{sF97}
Shmuel Friedland, \emph{The maximal orders of finite subgroups in {${\rm GL}\sb
  n({\bf Q})$}}, Proc. Amer. Math. Soc. \textbf{125} (1997), no.~12,
  3519--3526. \MR{98b:20064}

\bibitem{kImR93}
Kenneth Ireland and Michael Rosen, \emph{A classical introduction to modern
  number theory}, second ed., Graduate Texts in Mathematics, vol.~84,
  Springer-Verlag, New York, 1990. \MR{92e:11001}

\bibitem{mI76}
I.~Martin Isaacs, \emph{Character theory of finite groups}, Academic Press
  [Harcourt Brace Jovanovich Publishers], New York, 1976, Pure and Applied
  Mathematics, No. 69. \MR{57 \#417}

\bibitem{nJ85}
Nathan Jacobson, \emph{Basic algebra. {I}}, second ed., W. H. Freeman and
  Company, New York, 1985. \MR{86d:00001}

\bibitem{nJ89}
\bysame, \emph{Basic algebra. {II}}, second ed., W. H. Freeman and Company, New
  York, 1989. \MR{90m:00007}

\bibitem{gJ96}
Gerald~J. Janusz, \emph{Algebraic number fields}, second ed., Graduate Studies
  in Mathematics, vol.~7, American Mathematical Society, Providence, RI, 1996.
  \MR{96j:11137}

\bibitem{cJ78}
Camille Jordan, \emph{M\'emoire sur les {\'e}quations diff{\'e}rentielles
  lin{\'e}aires {\`a} int{\'e}grale alg{\'e}brique}, J. Reine Angew. Math.
  \textbf{84} (1878), 89--215.

\bibitem{yK94}
Yonatan~R. Katznelson, \emph{On the orders of finite subgroups of {${\rm
  GL}(n,\mathbb{Z})$}}, Exposition. Math. \textbf{12} (1994), 453--457.

\bibitem{KMRT}
Max-Albert Knus, Alexander Merkurjev, Markus Rost, and Jean-Pierre Tignol,
  \emph{The book of involutions}, American Mathematical Society Colloquium
  Publications, vol.~44, American Mathematical Society, Providence, RI, 1998,
  With a preface in French by J.\ Tits. \MR{2000a:16031}

\bibitem{jKaP02}
James Kuzmanovich and Andrey Pavlichenkov, \emph{Finite groups of matrices
  whose entries are integers}, Amer. Math. Monthly \textbf{109} (2002), no.~2,
  173--186. \MR{2003c:20057}

\bibitem{sL56}
Serge Lang, \emph{Algebraic groups over finite fields}, Amer. J. Math.
  \textbf{78} (1956), 555--563. \MR{19,174a}

\bibitem{sL90}
\bysame, \emph{Cyclotomic fields {I} and {II}}, second ed., Graduate Texts in
  Mathematics, vol. 121, Springer-Verlag, New York, 1990, With an appendix by
  Karl Rubin. \MR{91c:11001}

\bibitem{sL02}
\bysame, \emph{Algebra}, third ed., Graduate Texts in Mathematics, vol. 211,
  Springer-Verlag, New York, 2002. \MR{2003e:00003}

\bibitem{mLrP98}
Michael Larsen and Richard Pink, \emph{Finite subgroups of algebraic groups},
  preprint, 61 pages, 1998.

\bibitem{mL05}
Martin Lorenz, \emph{Multiplicative invariant theory}, Encyclopaedia of
  Mathematical Sciences, vol. 135, Springer-Verlag, Berlin, 2005, Invariant
  Theory and Algebraic Transformation Groups, VI.

\bibitem{hM85}
Hermann Minkowski, \emph{Untersuchungen {\"u}ber quadratische {F}ormen.
  {B}estimmung der {A}nzahl verschiedener {F}ormen, welche ein gegebenes
  {G}enus enth{\"a}lt}, Acta Mathematica \textbf{7} (1885), 201--258.

\bibitem{hM87}
\bysame, \emph{Zur {T}heorie der positiven quadratische {F}ormen}, J. reine
  angew. Math. \textbf{101} (1887), 196--202.

\bibitem{gNwP95}
Gabriele Nebe and Wilhelm Plesken, \emph{Finite rational matrix groups}, Mem.
  Amer. Math. Soc. \textbf{116} (1995), no.~556, viii+144. \MR{95k:20081}

\bibitem{dP71}
Donald~S. Passman, \emph{Infinite group rings}, Marcel Dekker Inc., New York,
  1971, Pure and Applied Mathematics, 6. \MR{47 \#3500}

\bibitem{dP77}
\bysame, \emph{The algebraic structure of group rings}, Wiley-Interscience
  [John Wiley \& Sons], New York, 1977, Pure and Applied Mathematics.
  \MR{81d:16001}

\bibitem{wP91}
Wilhelm Plesken, \emph{Some applications of representation theory},
  Representation theory of finite groups and finite-dimensional algebras
  (Bielefeld, 1991), Progr. Math., vol.~95, Birkh\"auser, Basel, 1991,
  pp.~477--496. \MR{92k:20019}

\bibitem{jjRot95}
Joseph~J. Rotman, \emph{An introduction to the theory of groups}, fourth ed.,
  Graduate Texts in Mathematics, vol. 148, Springer-Verlag, New York, 1995.
  \MR{95m:20001}

\bibitem{iS05}
Issai Schur, \emph{{\"U}ber eine {K}lasse von endlichen {G}ruppen linearer
  {S}ubstitutionen}, Sitzungsber. Preuss. Akad. Wiss. (1905), 77--91.

\bibitem{jpS73}
Jean-Pierre Serre, \emph{A course in arithmetic}, Springer-Verlag, New York,
  1973, Translated from the French, Graduate Texts in Mathematics, No. 7.
  \MR{49 \#8956}

\bibitem{jpS79}
\bysame, \emph{Local fields}, Graduate Texts in Mathematics, vol.~67,
  Springer-Verlag, New York, 1979, Translated from the French by Marvin Jay
  Greenberg. \MR{82e:12016}

\bibitem{OEIS}
Neil J.~A. Sloane, \emph{The on-line encyclopedia of integer sequences}, 2005,
  published electronically at
  \verb+http://www.research.att.com/~njas/sequences/+.

\bibitem{aS19}
Andreas Speiser, \emph{{Z}ahlentheoretische {S}{\"a}tze aus der
  {G}ruppentheorie}, Math. Zeitschrift \textbf{5} (1919), 1--6.

\bibitem{jW05}
Joseph H.~M. Wedderburn, \emph{A theorem on finite algebras}, Trans. Amer.
  Math. Soc. \textbf{6} (1905), 349--352.

\bibitem{bWxx}
Boris Weisfeiler, \emph{On the size and structure of finite linear groups},
  unfinished manuscript.

\bibitem{bW84}
\bysame, \emph{Post-classification version of {J}ordan's theorem on finite
  linear groups}, Proc. Nat. Acad. Sci. U.S.A. \textbf{81} (1984), no.~16,
  Phys. Sci., 5278--5279. \MR{85j:20041}

\end{thebibliography}


\end{document}